%% file: e3j.tex
\documentclass[12pt]{article}  
\usepackage{e-jc}

\usepackage{amsmath} 
\usepackage{latexsym}

\usepackage{url}
\usepackage{amssymb} 
\usepackage{exscale} 
\usepackage{graphicx} 
\usepackage{makeidx}

\usepackage[small,bf,hang]{caption} 

\usepackage{algorithm}
\usepackage{algorithmic}

\newtheorem{thm}{Theorem}

\newtheorem{lem}[thm]{Lemma}

\newenvironment{proof}[1][Proof]
{\par\noindent{\bf #1.} }{\hspace*{\fill}\nolinebreak{$\Box$}\bigskip\par}

\graphicspath{{figures/}}

\usepackage{color}

\title{\bf
The Ramsey Number $R(3,K_{10}-e)$\\
and Computational Bounds for $R(3,G)$
}

\author{
Jan Goedgebeur\footnote{
Supported by a Ph.D. grant from the Research Foundation of Flanders (FWO).}\\
\small Department of Applied Mathematics and Computer Science\small \\[-0.8ex]
\small Ghent University, B-9000 Ghent, Belgium\\[-0.8ex]
\small \texttt{jan.goedgebeur@gmail.com}\\
\\
Stanis\l aw P. Radziszowski\footnote{
Supported in part by the Polish National Science Centre grant 2011/02/A/ST6/00201.} \\
\small Department of Computer Science \\[-0.8ex]
\small Rochester Institute of Technology, Rochester, NY  14623, USA \\[-0.8ex]
\small \texttt{spr@cs.rit.edu}\\\
}

\date{\today\\
\dateline{September 1, 2013}{October 15, 2013}}

\begin{document}
\maketitle
\begin{abstract}
Using computer algorithms we establish that
the Ramsey number $R(3,K_{10}-e)$ is equal to 37, which
solves the smallest open case for Ramsey numbers of this type.
We also obtain new upper bounds for the cases of $R(3,K_k-e)$
for $11 \le k \le 16$, and show by construction a new lower
bound $55 \le R(3,K_{13}-e)$.

The new upper bounds on $R(3,K_k-e)$ are
obtained by using the values and lower bounds on $e(3,K_l-e,n)$
for $l \le k$, where
$e(3,K_k-e,n)$ is the minimum number of edges in any
triangle-free graph on $n$ vertices without $K_k-e$ in
the complement.
We complete the computation of the exact
values of $e(3,K_k-e,n)$ for all $n$ with $k \leq 10$ and for
$n \leq 34$ with $k = 11$, and establish many new lower
bounds on $e(3,K_k-e,n)$ for higher values of $k$.

Using the maximum triangle-free graph generation method,
we determine two other previously unknown Ramsey numbers,
namely $R(3,K_{10}-K_3-e)=31$ and $R(3,K_{10}-P_3-e)=31$.
For graphs $G$ on 10 vertices,
besides $G=K_{10}$, this leaves 6 open
cases of the form $R(3,G)$. The hardest among them
appears to be $G=K_{10}-2K_2$, for which we establish the
bounds $31 \le R(3,K_{10}-2K_2) \le 33$.

\end{abstract}


\input sect1.tex
\input sect2.tex
\input sect3.tex
\input sect4.tex
\input sect5.tex
\input ref.tex

\input appendix.tex

\end{document}

%% file: sect1.tex
\section{Notation and Preliminaries}

Notation, definitions and tools of this work are analogous to
those in our recent study of the classical two-color Ramsey
numbers $R(3,k)$ \cite{GoRa}, where $R(3,k)$ is
defined as the smallest $m$ such that no $m$-vertex
triangle-free graph with independence number
less than $k$ exists. Using coloring terminology,
such graphs can be
seen as 2-colorings of the edges of $K_m$, which have no
triangles in the first color and no
monochromatic $K_k$ in the second color. In this paper we study
mainly the case when $K_k-e$, the complete graph of order $k$
with one missing edge, is avoided in the second color.

\bigskip
Let $J_k$ denote the graph $K_k-e$, and hence the Ramsey
number $R(K_3,J_k)$ is the smallest $m$ such that every
triangle-free graph on $m$ vertices contains $J_k$ in the complement.
Similarly as in \cite{GoRa}, a graph $F$ will be called a
$(G,H;n,e)$-graph, if $|V(F)|=n, |E(F)|=e$, $F$ does
not contain $G$, and $\overline{F}$ (i.e.\ the complement of $F$)
does not contain $H$. By $\mathcal{R}(G,H;n,e)$ we denote
the set of all $(G,H;n,e)$-graphs.
We will often omit the parameter $e$, or both $e$ and $n$,
or give some range to either of these parameters, when referring
to special $(G,H;n,e)$-graphs or sets $\mathcal{R}(G,H;n,e)$.
If $G$ or $H$ is the complete graph $K_k$, we will often write
$k$ instead of $G$. For example, a $(3,J_k;n)$-graph $F$
is a $(K_3,K_k-e;n,e)$-graph for $e=|E(F)|$.

\medskip
In the remainder of this paper we will study only
triangle-free graphs, and mainly $K_k$ or $J_k$ will be avoided
in the complement. Note that for any $G \in \mathcal{R}(3,J_k)$
we have $\Delta(G)<k$, since all neighborhoods of vertices in $G$
are independent sets.
$e(3,J_k,n)$ $(=e(K_3,J_k,n))$ is defined as the
smallest number of edges in any $(3,J_k;n)$-graph.
The sum of the degrees of all neighbors of a vertex $v$ in $G$
will be denoted by $Z_G(v)$.
Similarly as in \cite{MPR,Ra1,ZoRa},
one can easily generalize the tools used in analysis
of the classical case $R(3,k)$ \cite{GoRa,GY,GR,RK1,RK2},
as described in the sequel.

\bigskip
Let $G$ be a $(3,J_k;n,e)$-graph.
For any vertex $v \in V(G)$, we will denote by $G_v$
the graph induced in $G$ by the set
$V(G) \setminus (N_G(v) \cup \{v\})$.
Note that if $d=\deg_G(v)$, then
$G_v$ is a $(3,J_{k-1};n-d-1,e-Z_G(v))$-graph.
This also implies that
$$\gamma(v)=\gamma(v,k,G)=e-Z_G(v)-e(3,J_{k-1},n-d-1)\ge 0,\eqno{(1)}$$
where $\gamma(v)$ is the so called {\em deficiency}
of vertex $v$ (as in \cite{GY}).
Finally, the deficiency of the graph $G$ is defined as
$$\gamma(G)=\sum_{v \in V(G)}{\gamma(v,k,G)} \ge 0.\eqno{(2)}$$

\noindent
The condition that $\gamma(G)\ge 0$ is often sufficient
to derive good lower bounds on $e(3,J_k,n)$, though a stronger
condition that all summands $\gamma(v,k,G)$
of (2) are non-negative sometimes implies better bounds.
It is easy to compute $\gamma(G)$ just from the degree
sequence of $G$ \cite{GY,GR,MPR}.
If $n_i$ is the number of vertices of degree $i$
in a $(3,J_k;n,e)$-graph $G$, then

\eject
$$\gamma(G)=ne-
\sum_{i}{n_i\big( i^2+e(3,J_{k-1},n-i-1)\big)} \ge 0,\eqno{(3)}$$
\noindent
where $n=\sum_{i=0}^{k-1}{n_i}$ and $2e=\sum_{i=0}^{k-1}{i n_i}$.

\bigskip
We obtain a number of improvements on lower bounds for
$e(3,J_k,n)$ and upper bounds for $R(3,J_k)$, summarized
at the end of the next section. The main computational result
of this paper solves the
smallest open case for the Ramsey numbers of the type
$R(3,J_k)$, namely we establish that $R(3,J_{10})=37$ by
improving the previous upper bound
$R(3,J_{10})\le 38$ \cite{MPR} by one.

Section~3 describes how
the algorithms of this work differed from those used
by us in the classical case~\cite{GoRa}, and how
we determined two other previously unknown Ramsey numbers,
namely $R(3,K_{10}-K_3-e)=31$ and $R(3,K_{10}-P_3-e)=31$,
using the maximum triangle-free graph generation method. 
Section~4 presents progress on $e(3,J_k,n)$ and
$R(3,J_k)$ for $k \le 11$, and Section~5
for $k \ge 12$.

\bigskip

%% file: sect2.tex
\section{Summary of Prior and New Results}

In 1995, Kim \cite{Kim} obtained 
a breakthrough result using probabilistic methods 
by establishing the exact asymptotics
for the classical case, namely
$R(3,k) = \Theta(n^2/\log n)$. The asymptotic behaviour of
$R(3,J_k)$ is clearly the same, since
$K_{k-1} \subset J_k \subset K_k$. The monotonicity
of $e(3,G,n)$ and Ramsey numbers $R(3,G)$ implies
that for all $n$ and $k$ we have
$$e(3,k,n)=e(K_3,K_k,n) \le e(K_3,J_k,n) \le e(K_3,K_{k-1},n),\eqno{(4)}$$
$$R(3,k-1)=R(K_3,K_{k-1}) \le R(K_3,J_k) \le R(K_3,K_k).\eqno{(5)}$$

\bigskip
For the small cases of $R(3,J_k)$ much of the progress
was obtained by deriving and using good lower bounds
on $e(3,J_k,n)$.
Explicit formulas for $e(3,J_{k+2},n)$ are known for all
$n \le 13k/4-1$, and for $n=13k/4$ when $k=0 \bmod 4$, as follows:

\medskip
\begin{thm}[\cite{ZoRa,RK1}]
\label{theorem:comulative_small}
For all $n,k \ge 1$, for which $e(3,J_{k+2},n)$ is finite, we have
$$
e(3,J_{k+2},n) = 
\left \{
\begin{array}{ll}
0 & \textrm{  if  }\  n \le k+1,\\
n - k & \textrm{  if  }\  k+2 \le n \le 2k \textrm{  and }\  k \ge 1,\\
3n - 5k & \textrm{  if  }\  2k < n \le 5k/2 \textrm{  and }\  k \ge 3,\\
5n - 10k & \textrm{  if  }\  5k/2 < n \le 3k \textrm{  and }\  k \ge 6,\\
6n - 13k & \textrm{  if  }\  3k < n \le 13k/4-1 \textrm{  and }\  k \ge 6.
\end{array}
\right. \eqno{(6)}
$$

\noindent
Furthermore, $e(3,J_{k+2},n) = 6n-13k$ for $k=4t$ and $n=13t$,
and the inequality $e(3,J_{k+2},n) \ge 6n-13k$ holds for all
$n$ and $k\ge 6$. All critical graphs have been
characterized whenever the equality in the theorem holds
for $n \le 3k$.
\end{thm}

\bigskip
\noindent
Our main focus in this direction is to obtain new exact
values and bounds on $e(3,J_{k+2},n)$ for $n \ge 13k/4$.
This in turn will permit us to prove the new upper
bounds on $R(3,J_k)$, for $10 \le k \le 16$.

\medskip
The general method we use is first to compute, if feasible,
the exact value of $e(3,J_k,n)$, or to derive a lower bound using
a combination of equalities (3) and (6), and computations. Better lower
bounds on $e(3,J_{k-1},m), m<n,$ often lead to better lower
bounds on $e(3,J_k,n)$. If we show that $e(3,J_k,n)=\infty$,
then we obtain an upper bound $R(3,J_k)\le n$.

Full enumeration of the sets $\mathcal{R}(3,J_k)$ for $k\le 6$
was completed in \cite{Ra1}, all such graphs for $k=7$ were
uploaded by Fidytek at a website~\cite{Fid},
and they were confirmed in this work.
Radziszowski computed the values of $e(3,J_7,n)$ and $e(3,J_8,n)$ in~\cite{Ra1}.
Some of the values and
bounds for $k=9$ and $k=10$, beyond those given by Theorem 1,
were obtained by McKay, Piwakowski and Radziszowski in~\cite{MPR}.
In this paper we complete this census for all cases of
$n$ with $k \le 10$, and give new
lower bounds for some higher parameters.

A $(3,J_k;n)$-graph is called \textit{critical} for a Ramsey
number $R(3,J_k)$ if $n = R(3,J_k) - 1$.
In~\cite{MPR}, McKay et al. determined that there are at
least 6 \textit{critical} triangle Ramsey graphs for $J_9$.
Using the maximum triangle-free method (see Section 3), we find
one more such graph and thus determine
that there are exactly 7 critical
graphs for $R(3,J_9)$. They can be downloaded from the
\textit{House of Graphs}~\cite{HOG} by searching for the
keywords ``critical ramsey graph for R(3,K9-e)''.

\medskip
There is an obvious similarity between Theorem 1 and
the results for $e(3,k,n)$ obtained in \cite{mtfg} as summarized
in Theorem 2 in \cite{GoRa}, though also note that there are some
differences. In particular, various cases are now restricted
to $k>c$. The graphs showing that these restrictions are
necessary are listed in \cite{ZoRa}.

\medskip
Our new results on $R(3,J_k)$ are marked in bold
in Table~\ref{table:new_bounds},
which presents the values and best known bounds on the Ramsey
numbers $R(3,J_k)$ and $R(3,K_k)$ for $k \le 16$.
The new upper bounds for $J_{10}$ and $J_{11}$ improve
on the bounds given in \cite{MPR} by 1 and 2, respectively.
Other upper bounds in bold are recorded for the first time.

\bigskip
\begin{table}[H]
\begin{center}
\begin{tabular}{|c|c|c||c|c|c|}
\hline
$k$&$R(3,J_k)$&$R(3,K_k)$&$k$&$R(3,J_k)$&$R(3,K_k)$\cr
\hline
3&\ \ 5&\ \ 6&10&{\bf 37}&40--42\cr
4&\ \ 7&\ \ 9&11&42--{\bf 45}&47--50\cr
5&11&14&12&47--{\bf 53}&52--59\cr
6&17&18&13&{\bf 55}--{\bf 62}&59--68\cr
7&21&23&14&59--{\bf 71}&66--77\cr
8&25&28&15&69--{\bf 80}&73--87\cr
9&31&36&16&73--{\bf 91}&82--98\cr
\hline
\end{tabular}

\caption{Ramsey numbers $R(3,J_k)$ and $R(3,K_k)$,
for $k \le 16$, $J_k=K_k-e$.}
\label{table:new_bounds}
\end{center}
\end{table}

The results $R(3,11)\ge 47$ \cite{Ex11} and
$R(3,16)\ge 82$ \cite{ExP} were recently obtained by Exoo.
Our recent work \cite{GoRa}, after 25 years of no progress, improved
the upper bound on $R(3,10)$ from 43 \cite{RK2} to 42, similarly as all
other upper bounds for $R(3,k)$ in the last column of Table 1.
The references for all other bounds and values,
and the previous bounds, are listed in \cite{SRN,GoRa}.

In a related cumulative work, Brinkmann, Goedgebeur
and Schlage-Puchta \cite{BGSP} completed the computation of all
Ramsey numbers of the form $R(3,G)$ for graphs $G$ on up
to 10 vertices, except 10 cases. The exceptions included
$K_{10}$, $J_{10}$, and 8 other graphs close to $K_{10}$.
The complements of these 10 graphs are depicted in Figure~1.
In fact,
the authors of~\cite{BGSP} showed in their article that the Ramsey number
for all of these remaining cases is at least 31.
In addition to $J_{10}$,
two of these cases are solved in this work,
namely $R(3,K_{10}-K_3-e)=R(3,K_{10}-P_3-e)=31$. Hence for graphs $G$
on 10 vertices, besides $G=K_{10}$, this leaves 6 other open
cases of the form $R(3,G)$. The hardest among them
appears to be $G=K_{10}-2K_2$, for which we establish the
bounds $31 \le R(3,K_{10}-2K_2) \le 33$.

\bigskip
\begin{figure}[h!t]
\centering
\includegraphics[width=.13\textwidth]{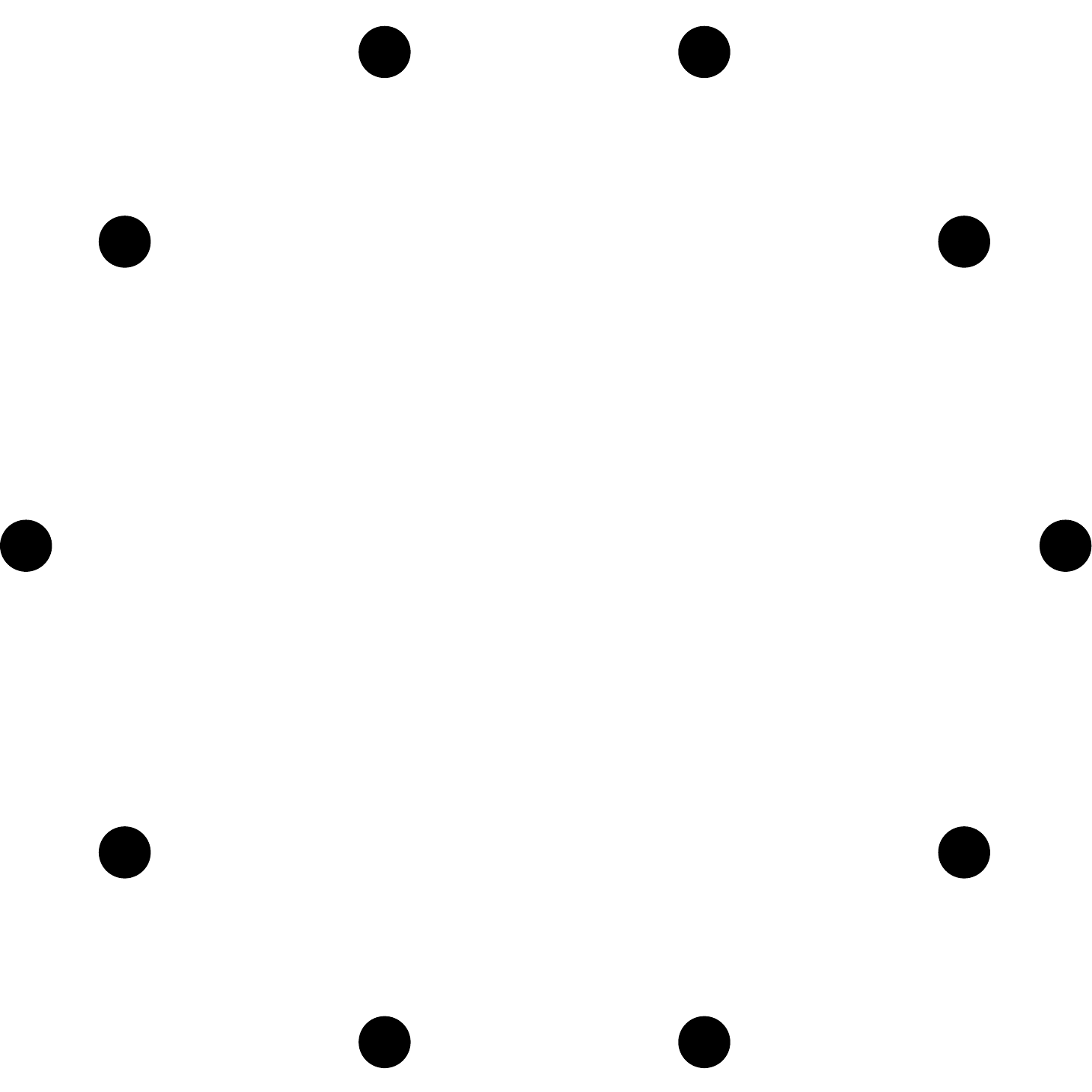}\ \ 
\includegraphics[width=.13\textwidth]{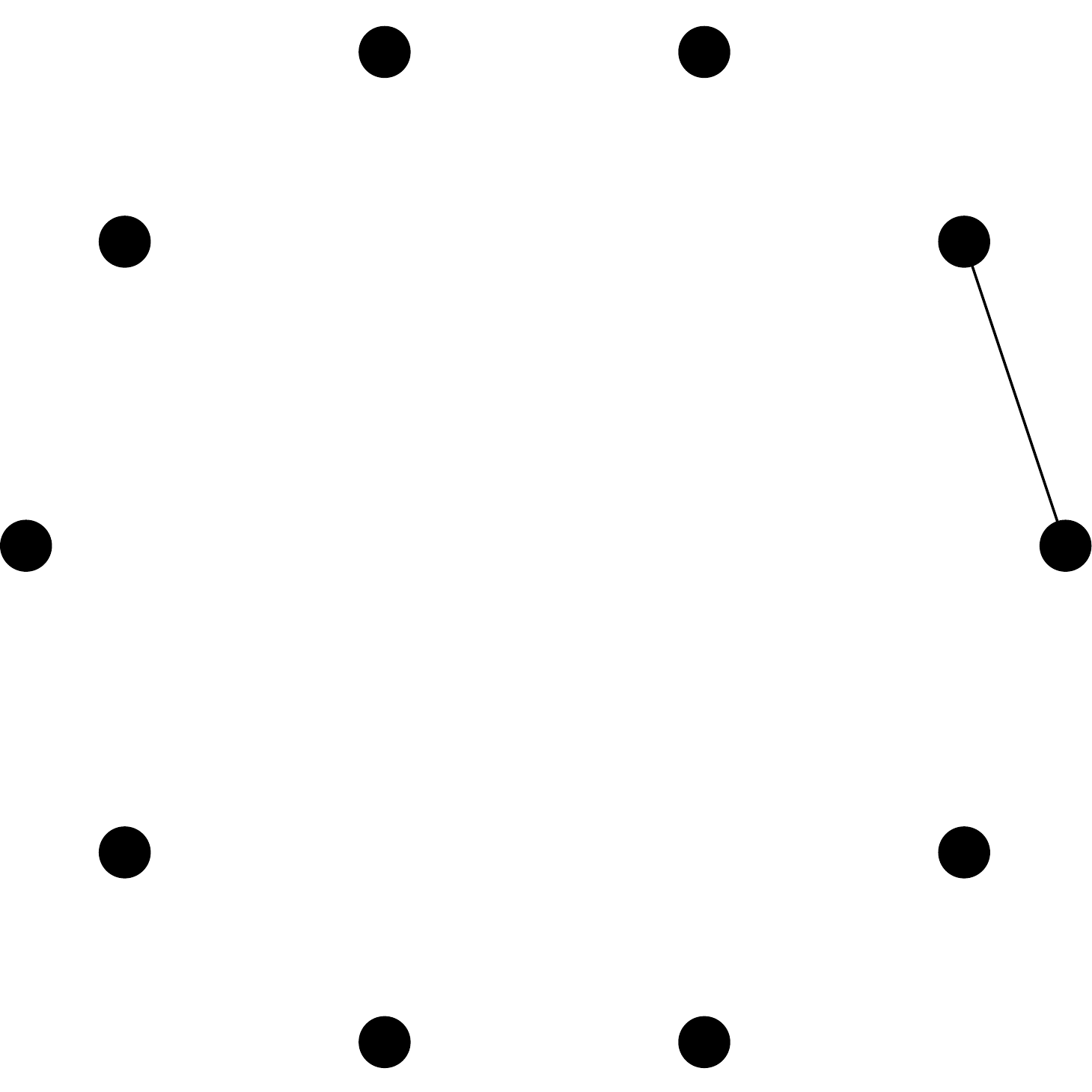}\ \ 
\includegraphics[width=.13\textwidth]{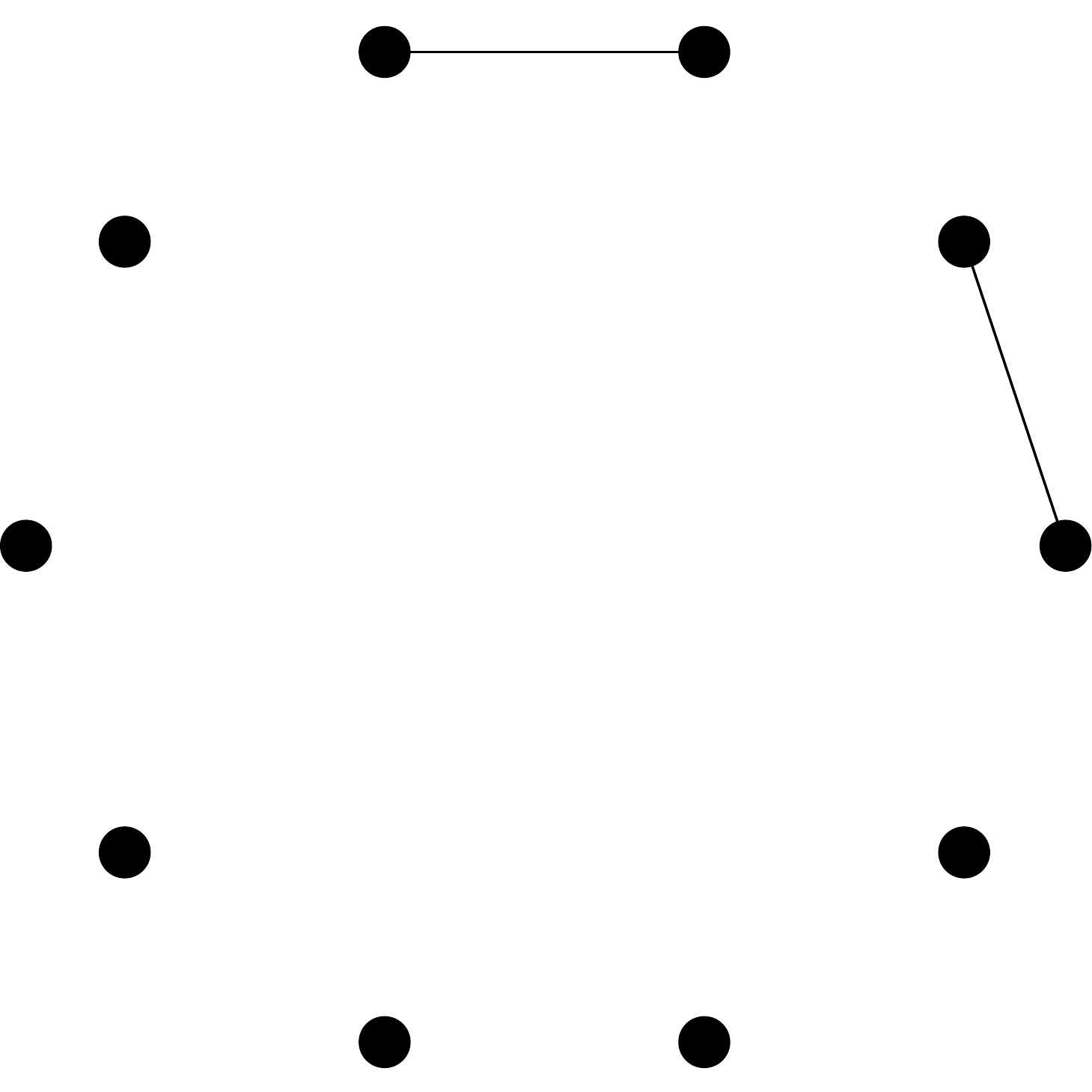}\ \ 
$\lfloor$
\includegraphics[width=.13\textwidth]{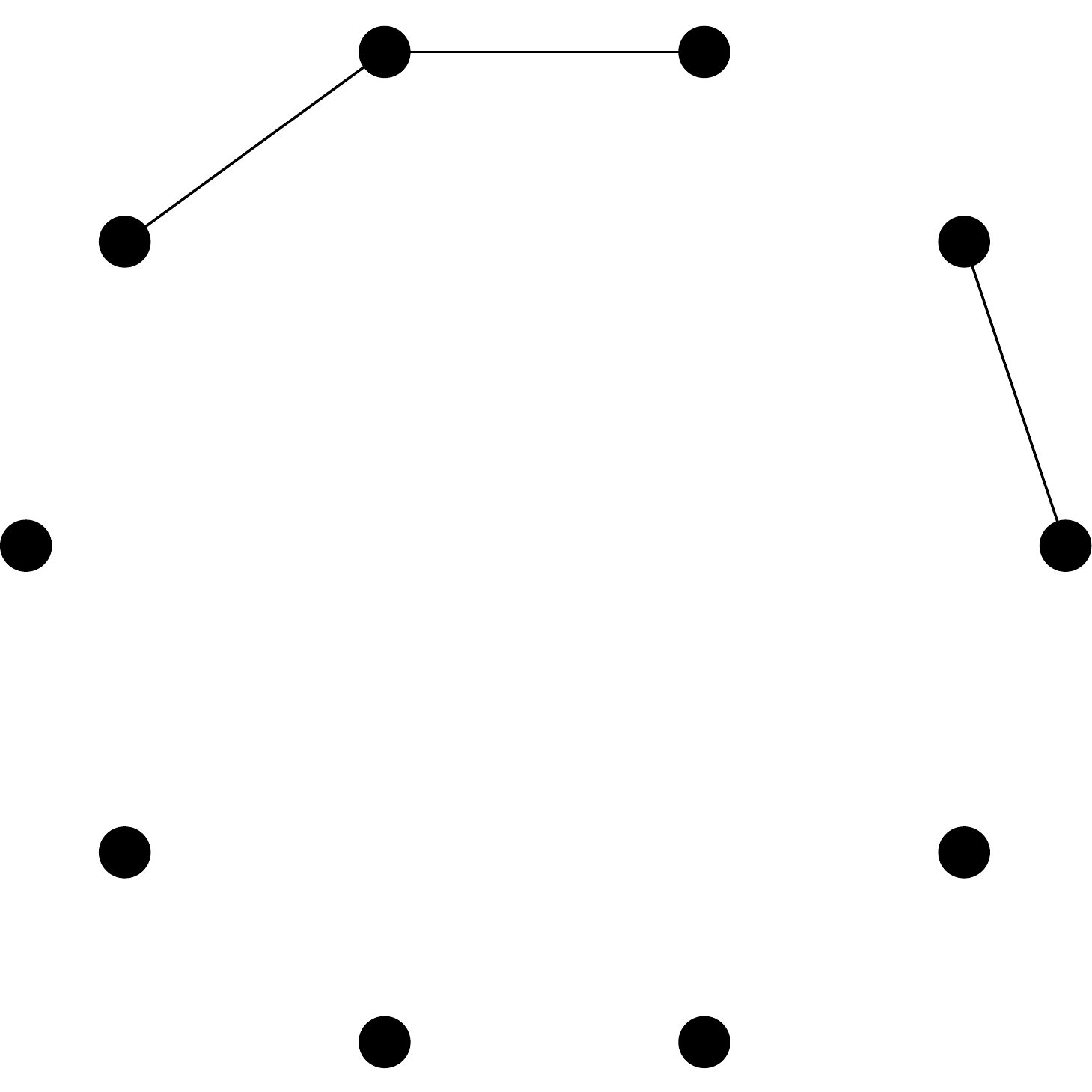}\ \ 
\includegraphics[width=.13\textwidth]{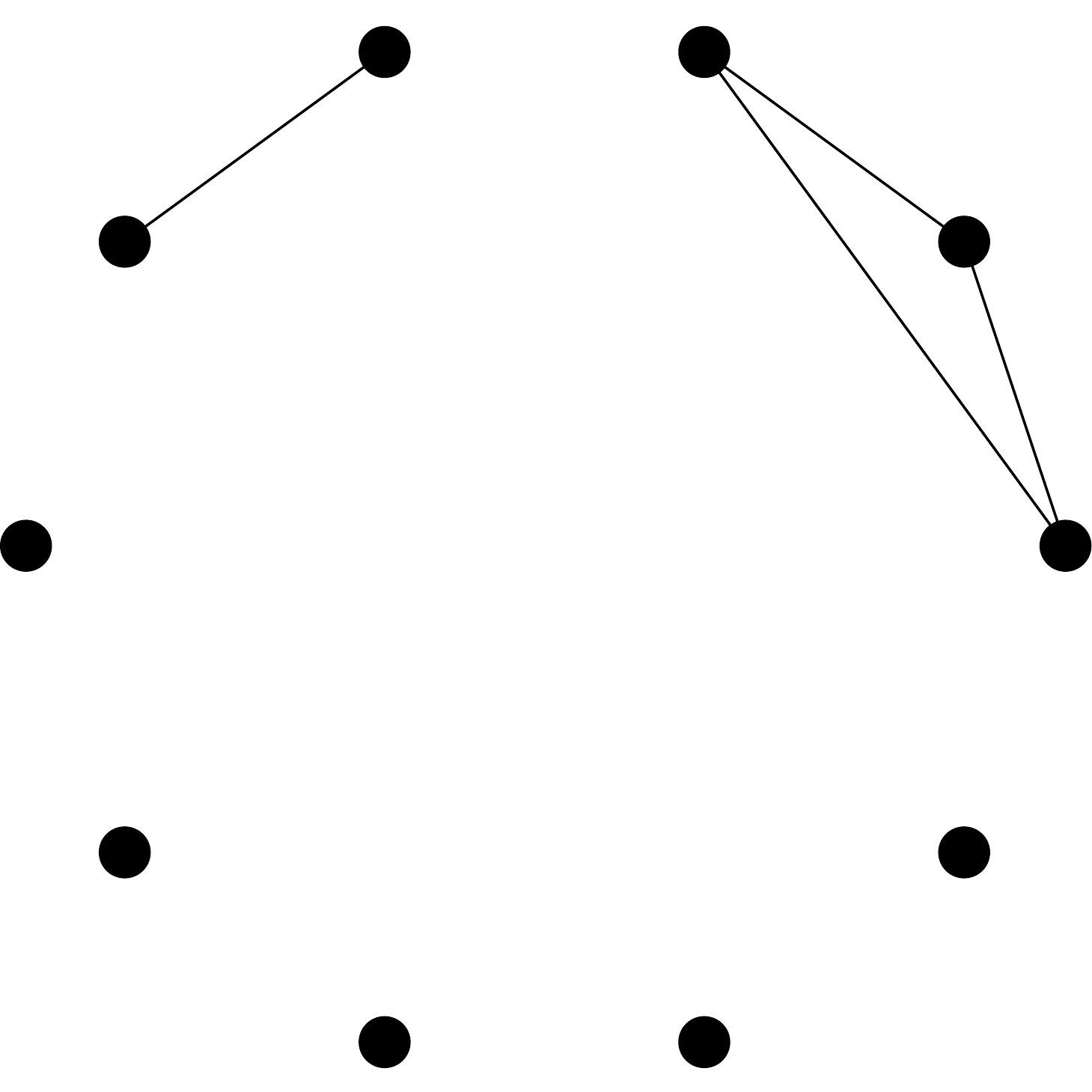}
$\rfloor$

\bigskip

\includegraphics[width=.13\textwidth]{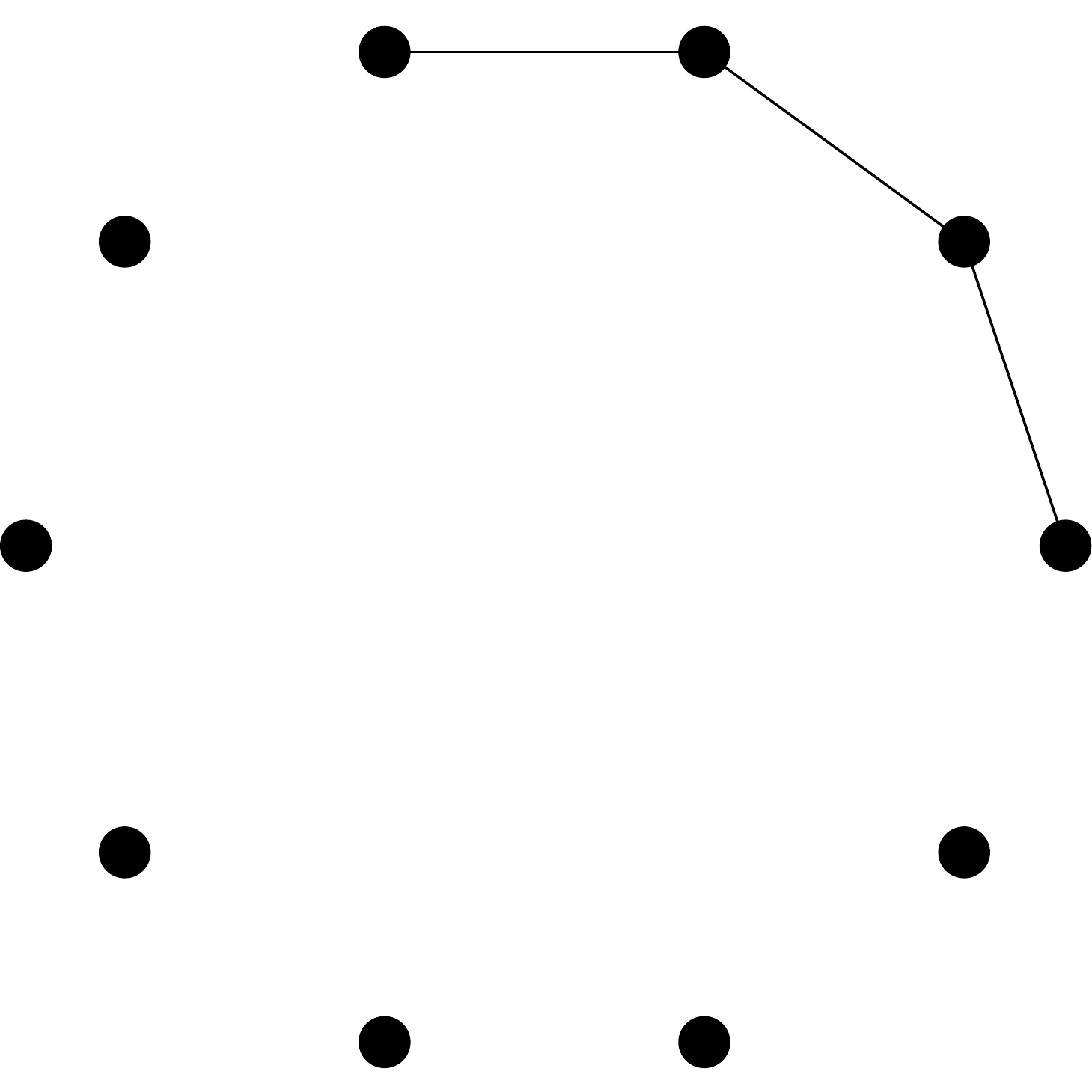}\ \ 
\includegraphics[width=.13\textwidth]{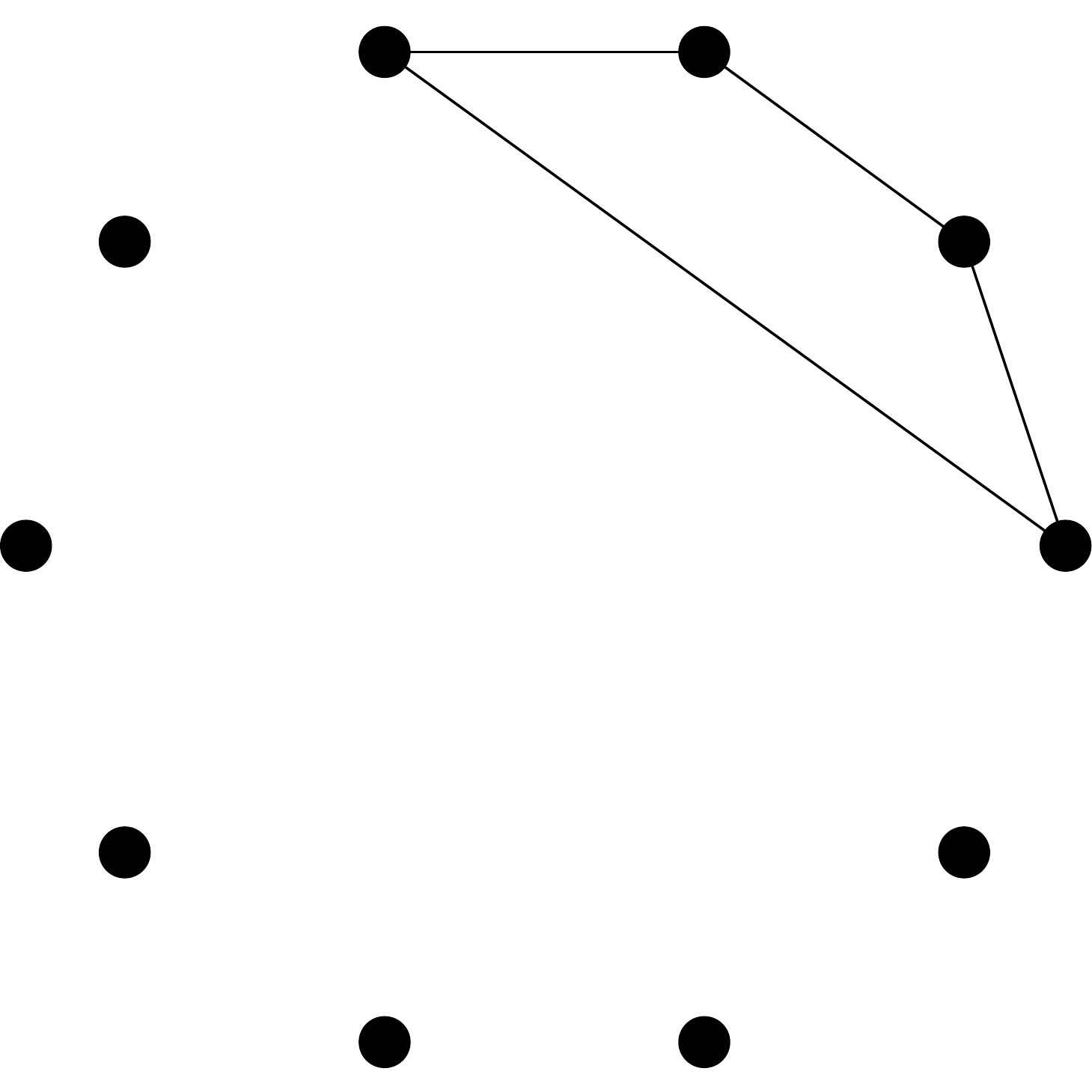}\ \ 
\includegraphics[width=.13\textwidth]{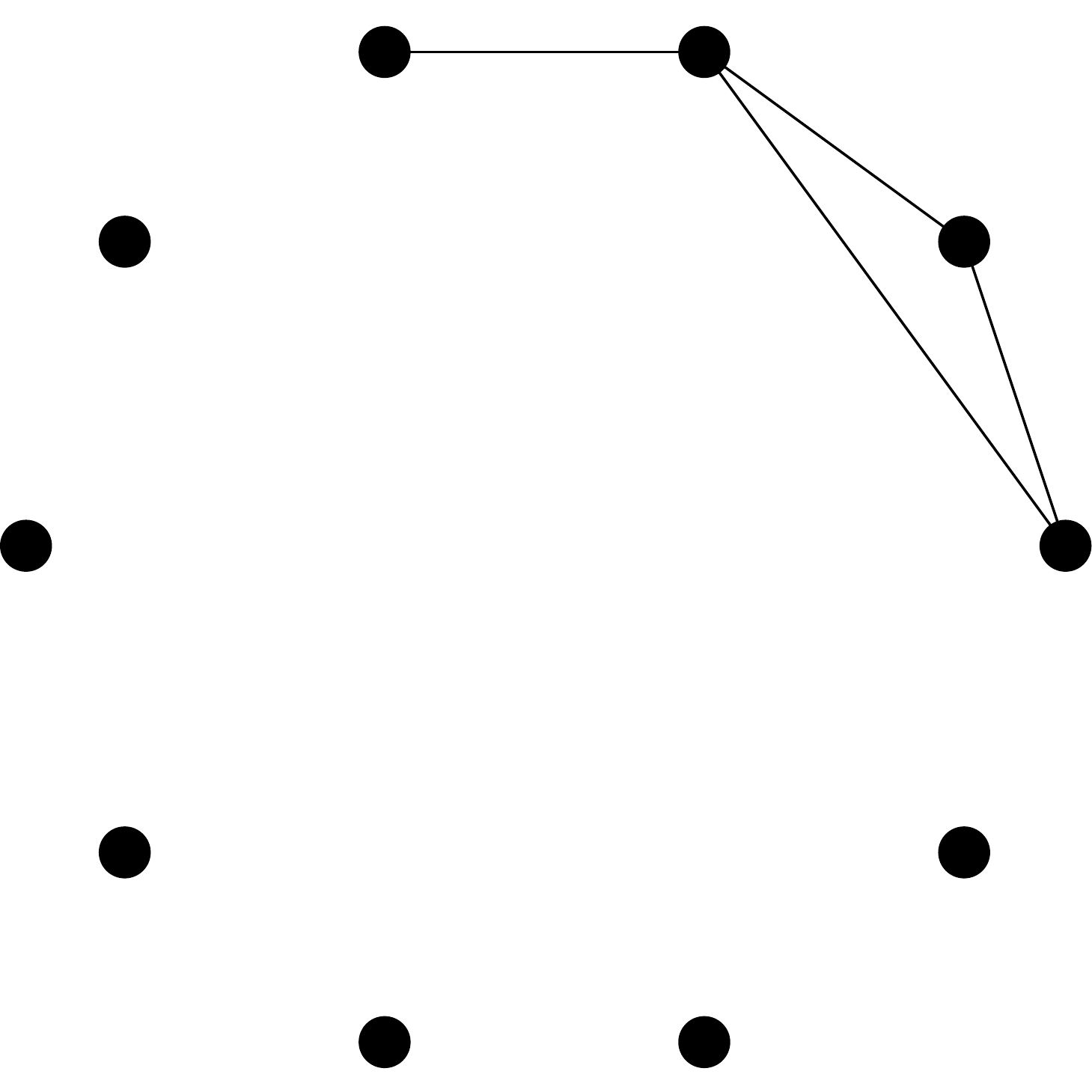}\ \ 
$\lfloor$
\includegraphics[width=.13\textwidth]{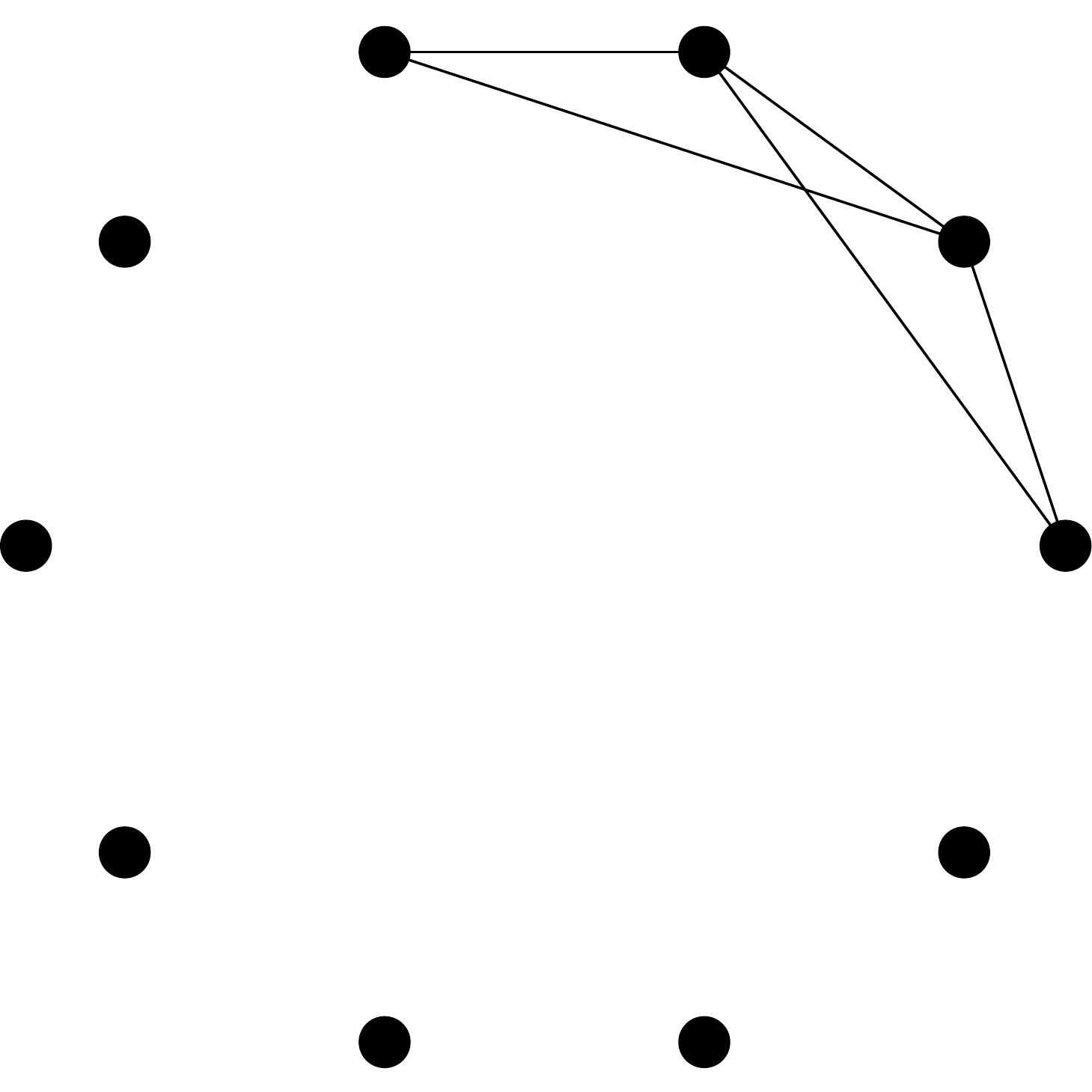}\ \ 
\includegraphics[width=.13\textwidth]{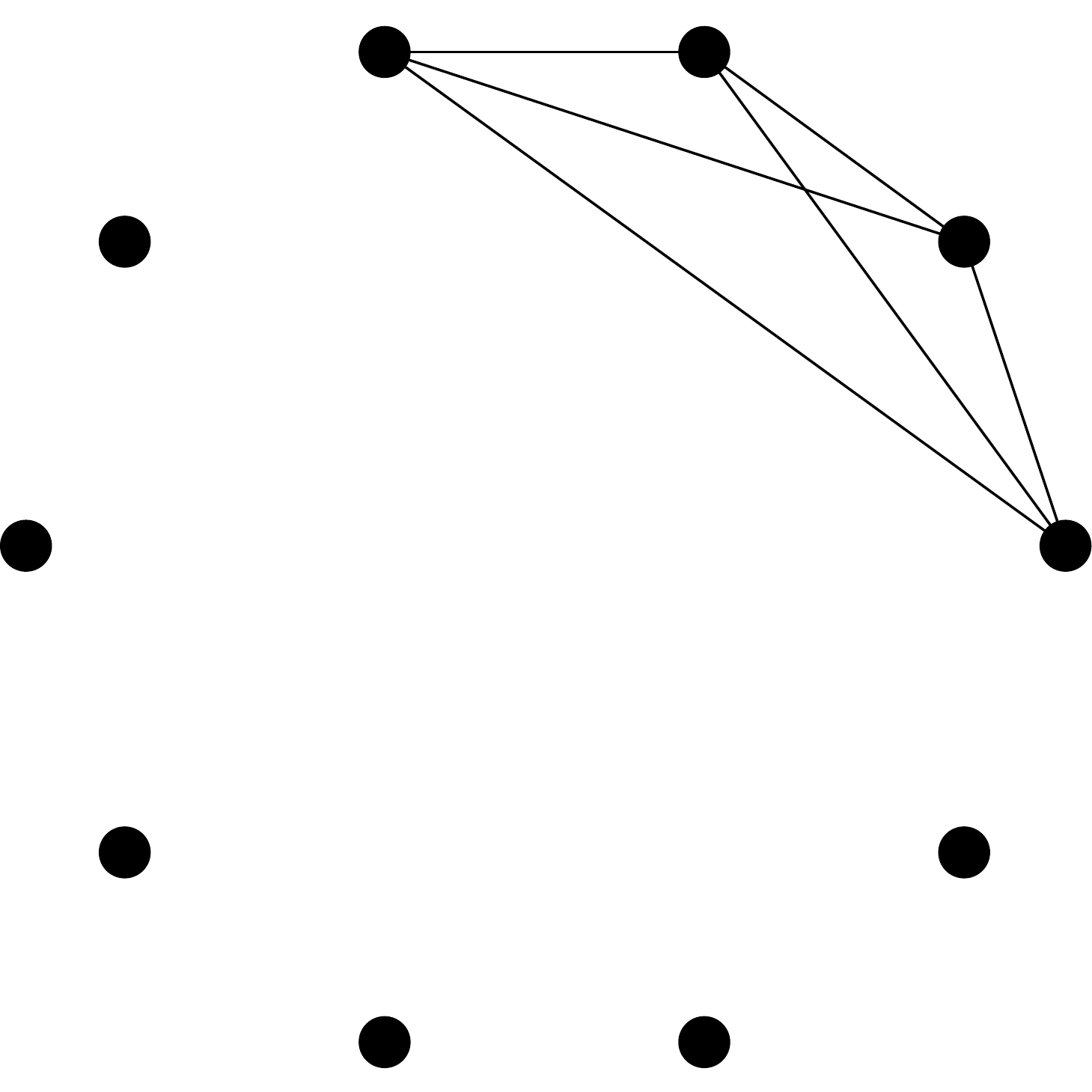}
$\rfloor$

\caption{The complements $\overline{G}$ of the 10 remaining graphs of order 10
which have $R(3,G) \ge 31$, for which Brinkmann et al.~\cite{BGSP} were unable to
determine the Ramsey number. Graphs which must have the same Ramsey number
are grouped by $\lfloor$ and $\rfloor$ (see~\cite{BGSP, Goed}).}
\label{fig:ramsey_10_unknown}
\end{figure}

\bigskip
We note
that in three cases of Table 1, namely for $k=12, 14$ and $16$,
the best known lower bounds for $R(3,J_k)$ are the same as for
$R(3,K_{k-1})$, and thus likely they can be improved.
Our new bounds on $R(3,J_k)$ are summarized in the following theorem. 

\begin{thm}
{\rm (a)}
$R(3,J_{10})=37$,
{\rm (b)}
$R(3,J_{11})\le 45$,
$R(3,J_{12})\le 53$,
$R(3,J_{13})\le 62$,
$R(3,J_{14})\le 71$,
$R(3,J_{15})\le 80$,
$R(3,J_{16})\le 91$, and
{\rm (c)}
$R(3,J_{13})\ge 55$.
\end{thm}

\begin{proof}
The lower bound $R(3,J_{10})\ge 37$ was established in \cite{MPR}.
The remaining sections describe our computational methods, and
intermediate values and bounds on $e(3,J_k,n)$. These imply
the upper bound for (a) in Section~4, and the bounds
for (b) in Section~5.
Result (c) follows from a circular $(3,J_{13};54)$-graph which we constructed.
It has arc distances $\{2, 3, 9, 16, 20, 24\}$. Note that all of these distances
have a nontrivial {\tt gcd} with 54, and thus this graph cannot be transformed
to an isomorphic circulant by a modular multiplier, so that 1 appears as one
of the distances.
\end{proof}

\bigskip
\noindent
{\bf Problem (Erd\H {o}s-S\'{o}s, 1980 \cite{Erd,CG})}
Let $\Delta_k = R(3,k)-R(3,k-1)$. Is it true that
$$\Delta_k \stackrel{k}{\rightarrow} \infty \  ? \ \ \ \ \
\Delta_k/k \stackrel{k}{\rightarrow} 0 \  ?$$

\medskip
Only easy bounds $3 \le \Delta_k \le k$ are known.
The results of this paper don't give any general improvement on the
bounds for $\Delta_k$, however we note that better understanding
of the behavior of $R(3,J_k)$ relative to $R(3,K_k)$ may lead
to such improvements since
$$\Delta_k = \big( R(3,K_k)-R(3,J_k) \big) + \big( R(3,J_k)-R(3,K_{k-1}) \big). \eqno{(7)}$$

\bigskip
The new results on $R(3,G)$ for some of the open cases
listed in Figure 1 are as follows.
\begin{thm}
{\rm (a)}
$31 \le R(3,K_{10}-2K_2) \le 33$,\\
{\rm (b)}
$R(3,K_{10}-K_3-e) = R(3,K_{10}-P_3-e) = 31$.
\end{thm}

\begin{proof}
The lower bound of 31 for each of the three cases was established in \cite{BGSP}.

The upper bound of 33 for (a) was obtained using essentially the same
method as the one used in the main case of this paper for $J_{10}$
(see Section~3 for more details). This improves over the trivial
bound of 37 implied by Theorem 2(a).
Table~\ref{table:graph_counts_K9-2e} in Appendix~1 contains information about data used
to derive the new bound: the values of $e(3,K_{9}-2K_2;n)$ and the counts
of corresponding graphs. Using only degree sequence analysis and the values
of $e(3,K_{9}-2K_2;n)$ one obtains the bound $R(3,K_{10}-2K_2) \le 35$.
Further computations using the neighbour gluing algorithm were required
to obtain the upper bound 33.

The computations applying the maximum triangle-free method from
\cite{BGSP} were enhanced
as described in Section 3, and gave the upper bounds needed for (b).
\end{proof}

Interestingly, contrary to our initial intuition, the case of
$K_{10}-2K_2$ appears to be significantly more difficult than $J_{10}$.
The computational effort which was required to prove $R(3,K_{10}-2K_2) \le 33$
was similar to the computational effort to prove $R(3,J_{10}) \le 37$,
but it looks like it is computationally infeasible to improve the upper bound
for $R(3,K_{10}-2K_2)$ any further by our current algorithms.
Our numerous attempts to improve the lower bound
failed, and consequently we conjecture that
$R(3,K_{10}-2K_2)=31$. If true, this would imply that for
each of the 6 remaining open cases of $G$ on 10 vertices (except $K_{10}$),
we have $R(3,G)=31$.

Finally, we would like to note that we performed exhaustive searches
for circulant graphs on up to 61 vertices in an attempt to improve
lower bounds for $R(3,k)$, $R(3,J_k)$, and for $R(3,G)$ for the remaining
graphs of order 10. If any of these lower bounds can still be
improved, it must be by using graphs which are not circulant.

%% file: sect3.tex
\bigskip
\section{Algorithms}
\label{section:algorithms}

Similarly as in~\cite{GoRa}, we use two independent techniques to
determine triangle Ramsey numbers: the maximum triangle-free method
and the neighborhood gluing extension method. These methods are
outlined in the following subsections.

\subsection*{Maximum Triangle-Free Method}

\subsubsection*{Generating maximal triangle-free graphs}
A maximal triangle-free graph (in short, an \textit{mtf} graph)
is a triangle-free graph such that the insertion of any new edge
forms a triangle. It is easy to see that there exists a $(3,J_k;n)$-graph
if and only if there is an mtf $(3,J_k;n)$-graph.
Brinkmann, Goedgebeur and Schlage-Puchta \cite{BGSP}
developed an algorithm to
exhaustively generate mtf graphs and mtf Ramsey graphs efficiently.
They implemented their algorithm in a program called
\textit{triangleramsey}~\cite{triangleramsey-site}. We refer the reader to~\cite{BGSP, Goed} for more details about the algorithm. Using this program they determined the
triangle Ramsey number $R(3,G)$ of nearly all graphs $G$ of order 10. The complements of the 10 graphs for which they were unable to determine the Ramsey number are
shown in Figure~\ref{fig:ramsey_10_unknown} in Section 2.

It is computationally infeasible to use \textit{triangleramsey} do
determine all mtf $(3,J_{10})$-graphs. However, we executed
\textit{triangleramsey} on a large computer cluster and were
able to determine all mtf Ramsey graphs for $K_{10}-P_3-e$ up
to 31 vertices (where $P_k$ is the path with $k$ vertices).
This is one of the remaining graphs whose Ramsey number could not be
determined by Brinkmann et al. The new computations took approximately
20 CPU years and the result is that there are 4 mtf Ramsey graphs
with 30 vertices for $K_{10}-P_3-e$
and no mtf Ramsey graphs with 31 vertices. Thus,
this proves that $R(3,K_{10}-P_3-e)=31$. 
By monotonicity of Ramsey numbers we have
$R(3,K_{10}-K_3-e) \le R(3,K_{10}-P_3-e)$, and thus the lower bound of 31 for both cases \cite{BGSP,Goed}
implies that $R(3,K_{10}-K_3-e)=R(3,K_{10}-P_3-e)=31$.

We also performed sample runs with \textit{triangleramsey} for the
other remaining graphs, but it looks like it will be computationally
infeasible to complete this task
by this method. E.g., we estimate that approximately 144 CPU years
would be required to determine $R(3,K_{10}-P_4)$ by running
\textit{triangleramsey}, and sample tests for
$R(3,K_{10}-K_4)$ indicate that this case will take
much longer than 200 CPU years.

\subsubsection*{Generating complete sets of triangle Ramsey graphs}
In order to determine $e(3,J_k,n)$ also non-maximal triangle-free
Ramsey graphs are required. Given all mtf $(3,J_k;n)$-graphs,
we can obtain all $(3,J_k;n)$-graphs by recursively removing edges
in all possible ways and testing if the obtained graphs are still
$(3,J_k;n)$-graphs. We used \textit{nauty}~\cite{McKay, nauty} to
make sure no isomorphic copies are output.
We generated, amongst others, the full sets
$\mathcal{R}(3,J_9;28), \mathcal{R}(3,J_9;29)$ and
$\mathcal{R}(3,J_9;30)$ (see Appendix~1 for detailed results)
using this method.

This mtf method is too slow for generating all $(3,J_{10};n)$-graphs
for $n$ which were needed in this work. Nevertheless, we used this
method to verify the correctness of our other programs for smaller
parameters. The results agreed in all cases in which more than one
method was used (see Appendix~2 for more details).

\subsection*{Neighborhood Gluing Extension Method}

The main method we used to improve upper bounds for
$R(3,J_k)$ is the neighborhood gluing extension method. In this
method our extension algorithm takes a  $(3,J_k;m)$-graph $H$ as
input and produces all $(3,J_{k+1};n,e)$-graphs $G$, often with
some specific restrictions on $n$ and $e$, such that for some vertex
$v \in V(G)$ graph $H$ is isomorphic to $G_v$. The program also gets
an expansion degree $d=n-m-1$ as input.
Thus it connects, or \textit{glues},
the $d$ neighbors of a vertex $v$ to $H$ in all possible ways. Note
that each neighbor of $v$ is glued to an independent set, otherwise
the extended graph would contain triangles. Similarly as in~\cite{GoRa},
various optimizations and bounding criteria are used to speed up
the algorithm. 

For example, suppose that we are aiming to construct $(3,J_{k+1})$-graphs.
Note that the complement of a graph $G$ contains
$J_k$ if and only if $G$ contains a spanning subgraph of $\overline{J_k}$
as an induced subgraph.
If two neighbors $u_1$ and $u_2$ of $v$ have already been connected
to independent sets $S_1$ and $S_2$ in $H$ and
$H[V(H) \setminus (S_1 \cup S_2)]$ contains a spanning subgraph of
$\overline{J_{k-1}}$ as induced subgraph, we can abort the recursion,
since this cannot yield any $(3,J_{k+1})$-graphs. 

There is, however, also one optimization which is specific to $J_k$.
Namely, we do not have to
connect the neighbors of $v$ to independent sets $S$ for which
$H[V(H) \setminus S]$ induces an independent set of order $k-1$,
since
otherwise this graph would contain $\overline{J_{k+1}}$ as induced
subgraph (an independent set of order $k-1$ together with
the disjoint edge $\{u,v\}$).

For more details about the general gluing algorithm,
we refer the reader to~\cite{GoRa,Goed}.

Most values and new bounds for $e(3,J_k,n)$, which are listed in
Section 4, were obtained by the gluing extension method.
In Appendix~2 we describe how we tested the correctness
of our implementation.

The strategy we used to determine if the parameters of the input graphs
to which our extender program was applied are sufficient, i.e. that
it is guaranteed that all $(3,J_{k+1};n,e)$-graphs are generated, is
the same as in~\cite{GoRa} and is outlined in the next subsection.

\subsection*{Degree Sequence Feasibility}

This method is based on the same principles as in the classical
case \cite{GoRa,RK2,Les}. Suppose we know the values
or lower bounds on $e(3,J_k,m)$ for some fixed $k$
and we wish to know all feasible degree sequences of $(3,J_{k+1};n,e)$-graphs.
We construct the system of integer constraints consisting of
$n=\sum_{i=0}^{k}{n_i}$, $2e=\sum_{i=0}^{k}{i n_i}$, and
the inequality (3).
If it has no solutions then we can conclude that no such
graphs exist.
Otherwise, we obtain solutions for $n_i$'s which include all
potential degree sequences.

%% file: sect4.tex
\bigskip
\section{Progress on Computing Small $e(3,J_k,n)$}
\label{section:computations_10}

\begin{table}[H]
\begin{center}
\begin{tabular}{|c|rrrrrrrrrrrrrr|}
\hline
vertices&\multicolumn{14}{c}{$k$}\vline\cr
$n$&3&4&5&6&7&8&9&10&11&12&13&14&15&16\cr
\hline
3&2&&&&&&&&&&&&&\cr
4&4&2&&&&&&&&&&&&\cr
5&$\infty$&4&2&&&&&&&&&&&\cr
6&&6&3&2&&&&&&&&&&\cr
7&&$\infty$&6&3&2&&&&&&&&&\cr
8&&&8&4&3&2&&&&&&&&\cr
9&&&12&7&4&3&2&&&&&&&\cr
10&&&15&10&5&4&3&2&&&&&&\cr
11&&&$\infty$&14&8&5&4&3&2&&&&&\cr
12&&&&18&11&6&5&4&3&2&&&&\cr
13&&&&24&15&9&6&5&4&3&2&&&\cr
14&&&&30&19&12&7&6&5&4&3&2&&\cr
15&&&&35&24&15&10&7&6&5&4&3&2&\cr
16&&&&40&30&20&13&8&7&6&5&4&3&2\cr
17&&&&$\infty$&37&25&16&11&8&7&6&5&4&3\cr
18&&&&&43&30&20&14&9&8&7&6&5&4\cr
19&&&&&54&37&25&17&12&9&8&7&6&5\cr
20&&&&&60&44&30&20&15&10&9&8&7&6\cr
21&&&&&$\infty$&51&35&25&18&13&10&9&8&7\cr
22&&&&&&59&{\bf 42}&30&21&16&11&10&9&8\cr
23&&&&&&70&{\bf 49}&35&25&19&14&11&10&9\cr
24&&&&&&80&{\bf 56}&40&30&22&17&12&11&10\cr
25&&&&&&$\infty$&{\bf 65}&46&35&25&20&15&12&11\cr
26&&&&&&&{\bf 73}&52&40&30&23&18&13&12\cr
27&&&&&&&{\bf 81}&{\bf 61}&45&35&26&21&16&13\cr
28&&&&&&&{\bf 95}&{\bf 68}&51&40&30&24&19&14\cr
29&&&&&&&{\bf 106}&{\bf 77}&{\bf 58}&45&35&27&22&17\cr
30&&&&&&&{\bf 117}&{\bf 86}&{\bf 66}&50&40&30&25&20\cr
31&&&&&&&$\infty$&{\bf 95}&{\bf 73}&56&45&35&28&23\cr
\hline
\end{tabular}

\caption{Exact values of $e(3,J_k,n)$,
for $3 \le k \le 16$, $3 \le n \le 31$.}
\label{table:exact_e_values_all}
\end{center}
\end{table}

\noindent
Most of the values of $e(3,J_k,n)$ collected in Table 2 are implied
by Theorem 1, others were obtained in \cite{Ra1,MPR,ZoRa}, and those
in bold are the result of this work.
The bottom-left blank part covers cases where the graphs with
corresponding parameters do not exist, while all entries in
the top-right blank area indicate 0 edges.
By Theorem 1 and \cite{mtfg}, for fixed $k$, $e(3,J_k,n)$ is equal
to $e(3,K_{k-1},n)$ for most small $n$, and from the data presented
in the following it looks like that this equality persists further
as $n$ grows. Only sporadic counterexamples to such behavior
for $n$ not much larger than $13k/4$ are known: seven such
cases are listed in \cite{ZoRa} for $k \le 7$,
and another one can be noted in Table 4 for $k=11, n=32$.
In other words, the second inequality of (4) seems to be 
much closer to equality than the first, and the opposite seems
to hold in (5). If true, we can expect that the first part
of $\Delta_k$ in (7) is significantly larger than the second part.

\bigskip
\noindent
{\bf Exact values of $e(3,J_9,n)$}

\medskip
\noindent
The values of $e(3,J_9,\le 21)$ are determined by Theorem 1.
The values of $e(3,J_9,n)$ for $22 \le n \le 30$ were obtained
by computations, mostly by the gluing extender algorithm which is outlined
in Section 3, and they are presented in Table 2.
These values improve over previously reported lower
bounds \cite{ZoRa,MPR}. We note that
$e(3,J_9,n)=e(3,K_8,n)$ for all $9 \le n \le 26$.

\bigskip
\noindent
{\bf Exact values of $e(3,J_{10},n)$}

\medskip
\noindent
The values of $e(3,J_{10},\le 26)$ are determined by Theorem 1.
The values for $27 \le n \le 37$ were obtained
by the gluing extender algorithm of Section 3,
and they are presented in Table 3.
These values improve over previously reported lower
bounds \cite{ZoRa,MPR}. We note that
$e(3,J_{10},n)=e(3,K_9,n)$ for all $10 \le n \le 35$
(see \cite{GoRa}).

\bigskip
\begin{table}[H]
\begin{center}
\begin{tabular}{c|c|l}
\hline
$n$&$e(3,J_{10},n)$&previous bound/comments\cr
\hline
26&\ \ 52&Theorem 1\cr
27&\ \ 61&58\cr
28&\ \ 68&65\cr
29&\ \ 77&72\cr
30&\ \ 86&81\cr
31&\ \ 95&90\cr
32&104&99\cr
33&118&110\cr
34&129&121\cr
35&140&133\cr
36&156&146, maximum 162\cr
37&$\infty$&hence $R(3,J_{10}) \le 37$, Theorem 2(a)\cr
\hline
\end{tabular}

\caption{Values of $e(3,J_{10},n)$, for $n\ge 26$.}
\label{table:bounds_e_10}
\end{center}
\end{table}

\noindent
{\bf Values and lower bounds on $e(3,J_{11},n)$}

\medskip
\noindent
Table 4 presents what we know about $e(3,J_{11},n)$ beyond Theorem 1,
which determines the values of $e(3,J_{11},\le 28)$.
The values and bounds for $29 \le n \le 41$ were obtained
by the gluing extender algorithm outlined in Section 3.
The lower bounds on $e(3,J_{11},\ge 42)$ are based on solving
integer constraints (2) and (3), using the exact values
of $e(3,J_{10},n)$ listed in Table 3. We note that
$e(3,J_{11},n)=e(3,K_{10},n)$ for all
$11 \le n \le 34$, except for $n=32$ (see \cite{GoRa}).
Four seemingly exceptional $(3,J_{11};32,80)$-graphs can
be obtained from the \textit{House of Graphs}~\cite{HOG}
by searching for the keywords ``exceptional minimal ramsey graph''.
One of them is the Wells graph, also called the Armanios-Wells
graph. It is a double cover of the complement of the Clebsch
graph. One of the other special graphs is formed by two disjoint
copies of the Clebsch graph itself. 
This works, since the Clebsch graph
is the unique $(3,J_6;16,40)$-graph \cite{Ra1} (denoted $G_4$
in the latter paper).

\bigskip
\begin{table}[H]
\begin{center}
\begin{tabular}{c|c|l}
\hline
$n$&$e(3,J_{11},n)\ge$&comments\cr
\hline
28&51&exact, Theorem 1\cr
29&58&exact\cr
30&66&exact\cr
31&73&exact\cr
32&80&exact, $e(3,10,32)=81$\cr
33&90&exact\cr
34&99&exact\cr
35&107&extender\cr
36&117&extender\cr
37&128&extender\cr
38&139&extender\cr
39&151&extender\cr
40&161&extender\cr
41&172&extender\cr
42&185&$e(3,10,42)=\infty$\cr
43&201&\cr
44&217&maximum 220\cr
45&$\infty$&hence $R(3,J_{11}) \le 45$, Theorem 2(b)\cr
\hline
\end{tabular}

\caption{Values and lower bounds on $e(3,J_{11},n)$, for $n\ge 28$.}
\label{table:bounds_e_11}
\end{center}
\end{table}

%% file: sect5.tex
\bigskip
\section{Progress on $e(3,J_k,n)$ and $R(3,J_k)$ for Higher $k$}

The results in Tables 3 and 4 required computations of
our gluing extender algorithm. We did not perform any such
computations in an attempt to improve the lower bounds on
$e(3,J_k,n)$ for $k \ge 12$, because such computations
would be hardly feasible. The results presented in
this section depend only on the degree sequence analysis
described in Section~3, using constraints (2), (3) and
the results for $k \le 11$ from the previous section.

\bigskip
\noindent
{\bf Lower bounds on $e(3,J_{12},n)$}

\noindent
Beyond the range of equality in Theorem 1 (for $n \ge 32$),
the lower bounds we obtained for
$e(3,J_{12},n)$ are the same as for
$e(3,K_{11},n)$, for all
$33 \le n \le 49$, except for $n=38$ (see \cite{GoRa}),
and they are presented in Table 5.
They were obtained by using constraints (2) and (3) as described
in Section~3. In one case, for $n=39$, we can
improve the lower bound by one as in the following lemma.

\medskip
\begin{lem}
$e(3,J_{12},39) \ge 117$.
\end{lem}
\begin{proof}
Suppose that $G$ is a $(3,J_{12};39,e)$-graph with $e \le 116$.
Using (2) and (3) with the bounds of Table 4 gives no solutions for
$e<116$ and two feasible degree sequences $n_i$ for $e=116$:
$n_4=1, n_6=38$ and $n_5=2, n_6=37$.
If $\deg(v)=4$ then $Z_G(v)=24$, and
if $\deg(v)=5$ then $Z_G(v) \ge 29$.
In both cases this contradicts inequality (1), and thus
we have $e(3,J_{12},39)\ge 117$.
\end{proof}

\bigskip
\bigskip
\begin{table}[H]
\begin{center}
\begin{tabular}{c|c|l}
\hline
$n$&$e(3,J_{12},n)\ge$&comments\cr
\hline
31&56&exact, Theorem 1\cr
32&62&Theorem 1\cr
33&68&Theorem 1\cr
34&75&\cr
35&83&\cr
36&92&\cr
37&100&\cr
38&108&$e(3,11,32)\ge 109$\cr
39&117&improvement by Lemma 4\cr
40&128&\cr
41&138&\cr
42&149&\cr
43&159&\cr
44&170&\cr
45&182&\cr
46&195&\cr
47&209&\cr
48&222&unique solution $n_7=36, n_8=12$\cr
49&237&\cr
50&252&$e(3,11,50)=\infty$\cr
51&266&\cr
52&280&maximum 286\cr
53&$\infty$&hence $R(3,J_{12}) \le 53$, Theorem 2(b)\cr
\hline
\end{tabular}

\caption{Lower bounds on $e(3,J_{12},n)$, for $n\ge 31$.}
\label{table:bounds_e_12}
\end{center}
\end{table}

\eject
\bigskip
\noindent
{\bf Lower bounds on $e(3,J_{13},n)$}

\noindent
Beyond the range of equality in Theorem 1 (for $n \ge 35$),
the lower bounds we obtained for $e(3,J_{13},n)$ are the same as for
$e(3,K_{12},n)$, for all
$35 \le n \le 58$, except for $n=45, 46$ (see \cite{GoRa}),
and they are presented in Table 6.
They were obtained by using constraints (2) and (3) as described
in Section 3.

\bigskip
\begin{table}[H]
\begin{center}
\begin{tabular}{c|c|l}
\hline
$n$&$e(3,J_{13},n)\ge$&comments\cr
\hline
34&61&exact, Theorem 1\cr
35&67&Theorem 1\cr
36&73&Theorem 1\cr
37&79&Theorem 1\cr
38&86&\cr
39&93&\cr
40&100&unique solution $n_5=40$\cr
41&109&unique solution $n_5=28, n_6=13$\cr
42&119&\cr
43&128&\cr
44&138&\cr
45&147&$e(3,12,45)\ge 148$\cr
46&157&$e(3,12,46)\ge 158$\cr
47&167&\cr
48&179&\cr
49&191&\cr
50&203&\cr
51&216&\cr
52&229&\cr
53&241&\cr
54&255&\cr
55&269&\cr
56&283&unique solution $n_{10}=50, n_{11}=6$\cr
57&299&\cr
58&316&\cr
59&333&$e(3,12,59)=\infty$\cr
60&350&maximum 360\cr
61&366&must be regular $n_{12}=61$\cr
62&$\infty$&hence $R(3,J_{13}) \le 62$, Theorem 2(b)\cr
\hline
\end{tabular}

\caption{Lower bounds on $e(3,J_{13},n)$, for $n\ge 34$.}
\label{table:bounds_e_13}
\end{center}
\end{table}

\eject
\bigskip
\noindent
{\bf Lower bounds on $e(3,J_{14},n)$}

\noindent
Beyond the range of equality in Theorem 1 (for $n \ge 40$),
the lower bounds we obtained for $e(3,J_{14},n)$ are the same as for
$e(3,K_{13},n)$, for all
$41 \le n \le 67$, except for $n=54, 55$ (see \cite{GoRa}),
and they are presented in Table 7.
They were obtained by using constraints (2) and (3) as described
in Section 3.

\bigskip
\begin{table}[H]
\begin{center}
\begin{tabular}{c|c|l}
\hline
$n$&$e(3,J_{14},n)\ge$&comments\cr
\hline
39&78&exact, Theorem 1\cr
40&84&Theorem 1\cr
41&91&\cr
42&97&\cr
43&104&\cr
44&112&\cr
45&120&\cr
46&128&\cr
47&136&\cr
48&146&\cr
49&157&\cr
50&167&\cr
51&177&\cr
52&189&\cr
53&200&\cr
54&210&$e(3,13,54)\ge 212$\cr
55&222&$e(3,13,55)\ge 223$\cr
56&234&\cr
57&247&\cr
58&260&\cr
59&275&\cr
60&289&\cr
61&303&\cr
62&319&\cr
63&334&\cr
64&350&\cr
65&365&\cr
66&381&\cr
67&398&\cr
68&416&$e(3,13,68)=\infty$\cr
69&434&\cr
70&451&maximum 455\cr
71&$\infty$&hence $R(3,J_{14}) \le 71$, Theorem 2(b)\cr
\hline
\end{tabular}

\caption{Lower bounds on $e(3,J_{14},n)$, for $n\ge 39$.}
\label{table:bounds_e_14}
\end{center}
\end{table}

\eject
\bigskip
\noindent
{\bf Lower bounds on $e(3,J_{15},n)$}

\noindent
The lower bounds we obtained for $e(3,J_{15},n)$
are the same as for $e(3,K_{14},n)$ for
all $71 \le n \le 76$ (see \cite{GoRa}),
and they are presented in Table 8.
They were obtained by using constraints (2) and (3) as described
in Section 3.

\bigskip
\begin{table}[H]
\begin{center}
\begin{tabular}{c|c|l}
\hline
$n$&$e(3,J_{15},n)\ge$&comments\cr
\hline
71&398&\cr
72&415&\cr
73&432&\cr
74&449&451 needed for $R(3,J_{16}) \le 90$\cr
75&468&\cr
76&486&473 sufficient for $R(3,J_{16}) \le 91$\cr
77&505&$e(3,14,77)=\infty$\cr
78&524&\cr
79&543&maximum 553\cr
80&$\infty$&hence $R(3,J_{15}) \le 80$, Theorem 2(b)\cr
\hline
\end{tabular}

\caption{Lower bounds on $e(3,J_{15},n)$, for $n\ge 71$.}
\label{table:bounds_e_15}
\end{center}
\end{table}

A 15-regular $(3,J_{16};90,675)$-graph $G$ is
feasible when for every vertex $v$ its $G_v$
is a $(3,J_{15};74,450)$-graph. Constraints (2) and (3)
have no feasible solution for $(3,J_{16};91)$-graphs,
and thus $R(3,J_{16}) \le 91$.

\bigskip
\subsection*{Acknowledgements}

We would like to thank Brendan McKay for his help in the search
for better lower bounds on $R(3,J_k)$.
Most computations for this work were carried out using
the Stevin Supercomputer Infrastructure at Ghent University.

%% file: ref.tex

%% file: appendix.tex
\bigskip
\bigskip
\section*{Appendix 1: Graph Counts}
\label{section:appendix1}

Tables~\ref{table:graph_counts_K7}--\ref{table:graph_counts_K11} below contain all known exact counts
of $(3,J_k;n,e)$-graphs for specified $n$,
for $k=7,8,9,10$ and $11$, respectively.
All graph counts were obtained by the algorithms described in
Section~3. Empty entries indicate 0.
In all cases, the maximum number of edges
is bounded by $\Delta(G)n/2 \le (k-1)n/2$.
All $(3,J_k;n, e(3,J_k,n))$-graphs, $k \le 10$,
which were constructed by our programs can be obtained from
the \textit{House of Graphs}~\cite{HOG} by searching for the
keywords ``minimal ramsey graph * Kk-e''.

All $(3,J_7)$-graphs were previously determined by Fidytek~\cite{Fid}.
We include the counts of $(3,J_7)$-graphs in
Table~\ref{table:graph_counts_K7} for completeness
and more uniform presentation,
since Fidytek provided statistics for $J_7$-free graphs
whose complement does not contain a $K_3$, while we list
triangle-free graphs whose complement does not contain a $J_7$.
The latter is similar to the
tables with data about $(3,J_k)$-graphs for $k \le 6$
in~\cite{Ra1}.

Table~\ref{table:graph_counts_K9-2e} contains all known exact counts
of $(3,K_9-2K_2;n,e)$-graphs. These graph counts were also obtained
by the algorithms described in Section~3. All edge minimal
$(3,K_9-2K_2;n, e(3,K_9-2K_2,n))$-graphs which were constructed
by our programs can be obtained from
the \textit{House of Graphs}~\cite{HOG} by searching for the
keywords ``minimal ramsey graph * K9-2K2''. 

\begingroup
\renewcommand{\baselinestretch}{0.95}

\begin{table}[H]
\begin{center}
{\scriptsize
\begin{tabular}{| c | ccccccccccccc |}
\hline 
edges & \multicolumn{13}{c}{number of vertices $n$} \vline\\
  $e$ & 8 & 9 & 10 & 11 & 12 & 13 & 14 & 15 & 16 & 17 & 18 & 19 & 20\\
\hline 
3 & 1  &    &    &    &    &    &    &    &    &    &    &    & \\
4  &  6  &  1  &    &    &    &    &    &    &    &    &    &    & \\
5  &  14  &  2  &  1  &    &    &    &    &    &    &    &    &    & \\
6  &  31  &  14  &  1  &    &    &    &    &    &    &    &    &    & \\
7  &  51  &  41  &  5  &    &    &    &    &    &    &    &    &    & \\
8  &  69  &  108  &  27  &  1  &    &    &    &    &    &    &    &    & \\
9  &  76  &  195  &  102  &  3  &    &    &    &    &    &    &    &    & \\
10  &  66  &  291  &  327  &  29  &    &    &    &    &    &    &    &    & \\
11  &  41  &  329  &  771  &  131  &  1  &    &    &    &    &    &    &    & \\
12  &  22  &  302  &  1355  &  643  &  8  &    &    &    &    &    &    &    & \\
13  &  9  &  204  &  1778  &  2158  &  47  &    &    &    &    &    &    &    & \\
14  &  3  &  117  &  1808  &  5239  &  398  &    &    &    &    &    &    &    & \\
15  &  2  &  53  &  1439  &  8961  &  2434  &  1  &    &    &    &    &    &    & \\
16  &  1  &  25  &  918  &  11450  &  9872  &  16  &    &    &    &    &    &    & \\
17  &    &  9  &  492  &  11072  &  26586  &  241  &    &    &    &    &    &    & \\
18  &    &  4  &  231  &  8505  &  49752  &  2665  &    &    &    &    &    &    & \\
19  &    &  1  &  99  &  5260  &  67226  &  16313  &  1  &    &    &    &    &    & \\
20  &    &  1  &  44  &  2794  &  68351  &  60891  &  13  &    &    &    &    &    & \\
21  &    &    &  19  &  1294  &  54124  &  145452  &  300  &    &    &    &    &    & \\
22  &    &    &  7  &  578  &  34707  &  238525  &  3997  &    &    &    &    &    & \\
23  &    &    &  3  &  233  &  18757  &  280341  &  28889  &    &    &    &    &    & \\
24  &    &    &  2  &  101  &  8976  &  247162  &  117123  &  2  &    &    &    &    & \\
25  &    &    &  1  &  41  &  3942  &  169011  &  291706  &  14  &    &    &    &    & \\
26  &    &    &    &  18  &  1669  &  93503  &  477533  &  305  &    &    &    &    & \\
27  &    &    &    &  6  &  693  &  43149  &  543408  &  4521  &    &    &    &    & \\
28  &    &    &    &  3  &  289  &  17392  &  451296  &  32828  &    &    &    &    & \\
29  &    &    &    &  1  &  115  &  6217  &  286635  &  121140  &    &    &    &    & \\
30  &    &    &    &  1  &  52  &  2073  &  146341  &  256923  &  3  &    &    &    & \\
31  &    &    &    &    &  21  &  626  &  63112  &  338238  &  22  &    &    &    & \\
32  &    &    &    &    &  10  &  190  &  24207  &  296128  &  361  &    &    &    & \\
33  &    &    &    &    &  4  &  50  &  8505  &  181637  &  3251  &    &    &    & \\
34  &    &    &    &    &  2  &  14  &  2841  &  83169  &  14968  &    &    &    & \\
35  &    &    &    &    &  1  &  3  &  884  &  30257  &  35296  &    &    &    & \\
36  &    &    &    &    &  1  &  1  &  275  &  9648  &  45855  &    &    &    & \\
37  &    &    &    &    &    &    &  75  &  2865  &  34944  &  1  &    &    & \\
38  &    &    &    &    &    &    &  22  &  883  &  16583  &  54  &    &    & \\
39  &    &    &    &    &    &    &  5  &  273  &  5269  &  349  &    &    & \\
40  &    &    &    &    &    &    &  2  &  94  &  1334  &  1070  &    &    & \\
41  &    &    &    &    &    &    &    &  32  &  350  &  1501  &    &    & \\
42  &    &    &    &    &    &    &    &  11  &  134  &  1174  &    &    & \\
43  &    &    &    &    &    &    &    &  4  &  50  &  522  &  2  &    & \\
44  &    &    &    &    &    &    &    &  1  &  25  &  147  &  8  &    & \\
45  &    &    &    &    &    &    &    &  1  &  8  &  26  &  38  &    & \\
46  &    &    &    &    &    &    &    &    &  4  &  6  &  61  &    & \\
47  &    &    &    &    &    &    &    &    &  1  &  1  &  58  &    & \\
48  &    &    &    &    &    &    &    &    &  1  &  1  &  36  &    & \\
49  &    &    &    &    &    &    &    &    &    &  1  &  17  &    & \\
50  &    &    &    &    &    &    &    &    &    &    &  4  &    & \\
51  &    &    &    &    &    &    &    &    &    &    &  1  &    & \\
52-53  &    &    &    &    &    &    &    &    &    &    &    &    & \\
54  &    &    &    &    &    &    &    &    &    &    &    &  1  & \\
55-59  &    &    &    &    &    &    &    &    &    &    &    &    & \\
60  &    &    &    &    &    &    &    &    &    &    &    &    &1\\
\hline
$|\mathcal{R}(3,J_7;n)|$ & 392 & 1697 & 9430 & 58522 & 348038 & 1323836 & 2447170 & 1358974 & 158459 & 4853 & 225 & 1 & 1 \\
\hline
\end{tabular}
}

\caption{Number of $(3,J_7;n,e)$-graphs, for $n \ge 8$.}

\label{table:graph_counts_K7}
\end{center}
\end{table}

\begin{table}[H]
\begin{center}
{\scriptsize
\begin{tabular}{| c | cccccccccc |}
\hline 
edges & \multicolumn{10}{c}{number of vertices $n$} \vline\\
  $e$ & 15 & 16 & 17 & 18 & 19 & 20 & 21 & 22 & 23 & 24\\
\hline 
15  &  1  &    &    &    &    &    &    &    &    &    \\
16  &  2  &    &    &    &    &    &    &    &    &    \\
17  &  18  &    &    &    &    &    &    &    &    &    \\
18  &  188  &    &    &    &    &    &    &    &    &    \\
19  &  ?  &    &    &    &    &    &    &    &    &    \\
20  &  ?  &  2  &    &    &    &    &    &    &    &    \\
21  &  ?  &  17  &    &    &    &    &    &    &    &    \\
22  &  ?  &  358  &    &    &    &    &    &    &    &    \\
23  &  ?  &  10659  &    &    &    &    &    &    &    &    \\
24  &  ?  &  ?  &    &    &    &    &    &    &    &    \\
25  &  ?  &  ?  &  2  &    &    &    &    &    &    &    \\
26  &  ?  &  ?  &  44  &    &    &    &    &    &    &    \\
27  &  ?  &  ?  &  2576  &    &    &    &    &    &    &    \\
28  &  ?  &  ?  &  117474  &    &    &    &    &    &    &    \\
29  &  ?  &  ?  &  ?  &    &    &    &    &    &    &    \\
30  &  ?  &  ?  &  ?  &  2  &    &    &    &    &    &    \\
31  &  ?  &  ?  &  ?  &  22  &    &    &    &    &    &    \\
32  &  ?  &  ?  &  ?  &  1175  &    &    &    &    &    &    \\
33  &  ?  &  ?  &  ?  &  79025  &    &    &    &    &    &    \\
34-36  &  ?  &  ?  &  ?  &  ?  &    &    &    &    &    &    \\
37  &  ?  &  ?  &  ?  &  ?  &  20  &    &    &    &    &    \\
38  &  ?  &  ?  &  ?  &  ?  &  2031  &    &    &    &    &    \\
39  &  ?  &  ?  &  ?  &  ?  &  130297  &    &    &    &    &    \\
40  &  ?  &  ?  &  ?  &  ?  &  3939009  &    &    &    &    &    \\
41-43  &  ?  &  ?  &  ?  &  ?  &  ?  &    &    &    &    &    \\
44  &  ?  &  ?  &  ?  &  ?  &  ?  &  169  &    &    &    &    \\
45  &  ?  &  ?  &  ?  &  ?  &  ?  &  8231  &    &    &    &    \\
46  &  ?  &  ?  &  ?  &  ?  &  ?  &  310400  &    &    &    &    \\
47  &  ?  &  ?  &  ?  &  ?  &  ?  &  5839714  &    &    &    &    \\
48-50  &  ?  &  ?  &  ?  &  ?  &  ?  &  ?  &    &    &    &    \\
51  &  ?  &  ?  &  ?  &  ?  &  ?  &  ?  &  7  &    &    &    \\
52  &  ?  &  ?  &  ?  &  ?  &  ?  &  ?  &  375  &    &    &    \\
53  &    &  ?  &  ?  &  ?  &  ?  &  ?  &  14141  &    &    &    \\
54  &    &  ?  &  ?  &  ?  &  ?  &  ?  &  255635  &    &    &    \\
55  &    &  ?  &  ?  &  ?  &  ?  &  ?  &  2262269  &    &    &    \\
56-58  &    &  ?  &  ?  &  ?  &  ?  &  ?  &  ?  &    &    &    \\
59  &    &    &  ?  &  ?  &  ?  &  ?  &  ?  &  2  &    &    \\
60  &    &    &    &  ?  &  ?  &  ?  &  ?  &  13  &    &    \\
61  &    &    &    &  ?  &  ?  &  ?  &  ?  &  162  &    &    \\
62  &    &    &    &  ?  &  ?  &  ?  &  ?  &  1630  &    &    \\
63  &    &    &    &  ?  &  ?  &  ?  &  ?  &  9101  &    &    \\
64  &    &    &    &    &  ?  &  ?  &  ?  &  26611  &    &    \\
65  &    &    &    &    &  ?  &  ?  &  ?  &  42700  &    &    \\
66  &    &    &    &    &  ?  &  ?  &  ?  &  41455  &    &    \\
67  &    &    &    &    &    &  ?  &  ?  &  26459  &    &    \\
68  &    &    &    &    &    &  ?  &  ?  &  11716  &    &    \\
69  &    &    &    &    &    &  ?  &  ?  &  3657  &    &    \\
70  &    &    &    &    &    &  ?  &  ?  &  957  &  1  &    \\
71  &    &    &    &    &    &    &  ?  &  208  &  2  &    \\
72  &    &    &    &    &    &    &  ?  &  42  &  8  &    \\
73  &    &    &    &    &    &    &  ?  &  10  &  6  &    \\
74  &    &    &    &    &    &    &    &  2  &  4  &    \\
75  &    &    &    &    &    &    &    &    &  1  &    \\
76  &    &    &    &    &    &    &    &    &  1  &    \\
77  &    &    &    &    &    &    &    &    &  13  &    \\
78-79  &    &    &    &    &    &    &    &    &    &    \\
80  &    &    &    &    &    &    &    &    &    &  1  \\
81-83  &    &    &    &    &    &    &    &    &    &  0  \\
84  &    &    &    &    &    &    &    &    &    &  8  \\
\hline
$|\mathcal{R}(3,J_8;n)|$  &  ?  &  ?  &  ?  &  ?  &  ?  &  ?  &  ?  &  164725  &  36  &  9  \\
\hline
\end{tabular}
}

\caption{Number of $(3,J_8;n,e)$-graphs, for $n \ge 15$.}

\label{table:graph_counts_K8}
\end{center}
\end{table}

\begin{table}[H]
\begin{center}
{\scriptsize
\begin{tabular}{| c | ccccccccccc |}
\hline 
edges & \multicolumn{11}{c}{number of vertices $n$} \vline\\
  $e$ & 20 & 21 & 22 & 23 & 24 & 25 & 26 & 27 & 28 & 29 & 30\\
\hline 
30  &  5  &    &    &    &    &    &    &    &    &    &    \\
31  &  64  &    &    &    &    &    &    &    &    &    &    \\
32  &  2073  &    &    &    &    &    &    &    &    &    &    \\
33-34  &  ?  &    &    &    &    &    &    &    &    &    &    \\
35  &  ?  &  1  &    &    &    &    &    &    &    &    &    \\
36  &  ?  &  20  &    &    &    &    &    &    &    &    &    \\
37  &  ?  &  951  &    &    &    &    &    &    &    &    &    \\
38  &  ?  &  39657  &    &    &    &    &    &    &    &    &    \\
39-41  &  ?  &  ?  &    &    &    &    &    &    &    &    &    \\
42  &  ?  &  ?  &  21  &    &    &    &    &    &    &    &    \\
43  &  ?  &  ?  &  1592  &    &    &    &    &    &    &    &    \\
44  &  ?  &  ?  &  86833  &    &    &    &    &    &    &    &    \\
45  &  ?  &  ?  &  3963053  &    &    &    &    &    &    &    &    \\
46-48  &  ?  &  ?  &  ?  &    &    &    &    &    &    &    &    \\
49  &  ?  &  ?  &  ?  &  103  &    &    &    &    &    &    &    \\
50  &  ?  &  ?  &  ?  &  9102  &    &    &    &    &    &    &    \\
51  &  ?  &  ?  &  ?  &  514099  &    &    &    &    &    &    &    \\
52-55  &  ?  &  ?  &  ?  &  ?  &    &    &    &    &    &    &    \\
56  &  ?  &  ?  &  ?  &  ?  &  54  &    &    &    &    &    &    \\
57  &  ?  &  ?  &  ?  &  ?  &  3639  &    &    &    &    &    &    \\
58  &  ?  &  ?  &  ?  &  ?  &  173608  &    &    &    &    &    &    \\
59-64  &  ?  &  ?  &  ?  &  ?  &  ?  &    &    &    &    &    &    \\
65  &  ?  &  ?  &  ?  &  ?  &  ?  &  547  &    &    &    &    &    \\
66  &  ?  &  ?  &  ?  &  ?  &  ?  &  48964  &    &    &    &    &    \\
67  &  ?  &  ?  &  ?  &  ?  &  ?  &  2538589  &    &    &    &    &    \\
68-72  &  ?  &  ?  &  ?  &  ?  &  ?  &  ?  &    &    &    &    &    \\
73  &  ?  &  ?  &  ?  &  ?  &  ?  &  ?  &  62  &    &    &    &    \\
74  &  ?  &  ?  &  ?  &  ?  &  ?  &  ?  &  1857  &    &    &    &    \\
75  &  ?  &  ?  &  ?  &  ?  &  ?  &  ?  &  36799  &    &    &    &    \\
76  &  ?  &  ?  &  ?  &  ?  &  ?  &  ?  &  755052  &    &    &    &    \\
77-80  &  ?  &  ?  &  ?  &  ?  &  ?  &  ?  &  ?  &    &    &    &    \\
81  &    &  ?  &  ?  &  ?  &  ?  &  ?  &  ?  &  4  &    &    &    \\
82  &    &  ?  &  ?  &  ?  &  ?  &  ?  &  ?  &  24  &    &    &    \\
83  &    &  ?  &  ?  &  ?  &  ?  &  ?  &  ?  &  197  &    &    &    \\
84  &    &  ?  &  ?  &  ?  &  ?  &  ?  &  ?  &  1126  &    &    &    \\
85  &    &    &  ?  &  ?  &  ?  &  ?  &  ?  &  6206  &    &    &    \\
86  &    &    &  ?  &  ?  &  ?  &  ?  &  ?  &  42468  &    &    &    \\
87  &    &    &  ?  &  ?  &  ?  &  ?  &  ?  &  384398  &    &    &    \\
88  &    &    &  ?  &  ?  &  ?  &  ?  &  ?  &  2843005  &    &    &    \\
89-94  &    &    &    &  ?  &  ?  &  ?  &  ?  &  ?  &    &    &    \\
95  &    &    &    &    &  ?  &  ?  &  ?  &  ?  &  1  &    &    \\
96  &    &    &    &    &  ?  &  ?  &  ?  &  ?  &  14  &    &    \\
97  &    &    &    &    &    &  ?  &  ?  &  ?  &  107  &    &    \\
98  &    &    &    &    &    &  ?  &  ?  &  ?  &  1062  &    &    \\
99  &    &    &    &    &    &  ?  &  ?  &  ?  &  5182  &    &    \\
100  &    &    &    &    &    &  ?  &  ?  &  ?  &  16588  &    &    \\
101  &    &    &    &    &    &    &  ?  &  ?  &  34077  &    &    \\
102  &    &    &    &    &    &    &  ?  &  ?  &  50241  &    &    \\
103  &    &    &    &    &    &    &  ?  &  ?  &  51686  &    &    \\
104  &    &    &    &    &    &    &  ?  &  ?  &  39702  &    &    \\
105  &    &    &    &    &    &    &    &  ?  &  21621  &    &    \\
106  &    &    &    &    &    &    &    &  ?  &  9379  &  1  &    \\
107  &    &    &    &    &    &    &    &  ?  &  2864  &  0  &    \\
108  &    &    &    &    &    &    &    &  ?  &  843  &  0  &    \\
109  &    &    &    &    &    &    &    &    &  158  &  2  &    \\
110  &    &    &    &    &    &    &    &    &  49  &  6  &    \\
111  &    &    &    &    &    &    &    &    &  7  &  9  &    \\
112  &    &    &    &    &    &    &    &    &  91  &  6  &    \\
113-115  &    &    &    &    &    &    &    &    &    &  0  &    \\
116  &    &    &    &    &    &    &    &    &    &  1  &    \\
117  &    &    &    &    &    &    &    &    &    &    &  1  \\
118  &    &    &    &    &    &    &    &    &    &    &  1  \\
119  &    &    &    &    &    &    &    &    &    &    &  1  \\
120  &    &    &    &    &    &    &    &    &    &    &  4  \\
\hline
$|\mathcal{R}(3,J_9;n)|$  &  ?  &  ?  &  ?  &  ?  &  ?  &  ?  &  ?  &  ?  &  233672  &  25  &  7  \\
\hline
\end{tabular}
}

\caption{Number of $(3,J_9;n,e)$-graphs, for $n \ge 20$.}

\label{table:graph_counts_K9}
\end{center}
\end{table}

\begin{table}[H]
\begin{center}
{\scriptsize
\begin{tabular}{| c | ccccccccccccc |}
\hline 
edges & \multicolumn{13}{c}{number of vertices $n$} \vline\\
  $e$ & 24 & 25 & 26 & 27 & 28 & 29 & 30 & 31 & 32 & 33 & 34 & 35 & 36\\
\hline 
40  &  2  &    &    &    &    &    &    &    &    &    &    &    &    \\
41  &  32  &    &    &    &    &    &    &    &    &    &    &    &    \\
42  &  2089  &    &    &    &    &    &    &    &    &    &    &    &    \\
43-45  &  ?  &    &    &    &    &    &    &    &    &    &    &    &    \\
46  &  ?  &  1  &    &    &    &    &    &    &    &    &    &    &    \\
47  &  ?  &  39  &    &    &    &    &    &    &    &    &    &    &    \\
48  &  ?  &  4113  &    &    &    &    &    &    &    &    &    &    &    \\
49-51  &  ?  &  ?  &    &    &    &    &    &    &    &    &    &    &    \\
52  &  ?  &  ?  &  1  &    &    &    &    &    &    &    &    &    &    \\
53  &  ?  &  ?  &  1  &    &    &    &    &    &    &    &    &    &    \\
54  &  ?  &  ?  &  444  &    &    &    &    &    &    &    &    &    &    \\
55  &  ?  &  ?  &  58550  &    &    &    &    &    &    &    &    &    &    \\
56-60  &  ?  &  ?  &  ?  &    &    &    &    &    &    &    &    &    &    \\
61  &  ?  &  ?  &  ?  &  700  &    &    &    &    &    &    &    &    &    \\
62  &  ?  &  ?  &  ?  &  95185  &    &    &    &    &    &    &    &    &    \\
63  &  ?  &  ?  &  ?  &  6531339  &    &    &    &    &    &    &    &    &    \\
64-67  &  ?  &  ?  &  ?  &  ?  &    &    &    &    &    &    &    &    &    \\
68  &  ?  &  ?  &  ?  &  ?  &  126  &    &    &    &    &    &    &    &    \\
69  &  ?  &  ?  &  ?  &  ?  &  17223  &    &    &    &    &    &    &    &    \\
70  &  ?  &  ?  &  ?  &  ?  &  1204171  &    &    &    &    &    &    &    &    \\
71-76  &  ?  &  ?  &  ?  &  ?  &  ?  &    &    &    &    &    &    &    &    \\
77  &  ?  &  ?  &  ?  &  ?  &  ?  &  1342  &    &    &    &    &    &    &    \\
78  &  ?  &  ?  &  ?  &  ?  &  ?  &  156982  &    &    &    &    &    &    &    \\
79-85  &  ?  &  ?  &  ?  &  ?  &  ?  &  ?  &    &    &    &    &    &    &    \\
86  &  ?  &  ?  &  ?  &  ?  &  ?  &  ?  &  1800  &    &    &    &    &    &    \\
87  &  ?  &  ?  &  ?  &  ?  &  ?  &  ?  &  147407  &    &    &    &    &    &    \\
88-94  &  ?  &  ?  &  ?  &  ?  &  ?  &  ?  &  ?  &    &    &    &    &    &    \\
95  &  ?  &  ?  &  ?  &  ?  &  ?  &  ?  &  ?  &  560  &    &    &    &    &    \\
96  &  ?  &  ?  &  ?  &  ?  &  ?  &  ?  &  ?  &  35154  &    &    &    &    &    \\
97-103  &  ?  &  ?  &  ?  &  ?  &  ?  &  ?  &  ?  &  ?  &    &    &    &    &    \\
104  &  ?  &  ?  &  ?  &  ?  &  ?  &  ?  &  ?  &  ?  &  39  &    &    &    &    \\
105  &  ?  &  ?  &  ?  &  ?  &  ?  &  ?  &  ?  &  ?  &  952  &    &    &    &    \\
106-117  &  ?  &  ?  &  ?  &  ?  &  ?  &  ?  &  ?  &  ?  &  ?  &    &    &    &    \\
118  &    &    &    &  ?  &  ?  &  ?  &  ?  &  ?  &  ?  &  $\ge 5$  &    &    &    \\
119  &    &    &    &  ?  &  ?  &  ?  &  ?  &  ?  &  ?  &  $\ge 86$  &    &    &    \\
120  &    &    &    &  ?  &  ?  &  ?  &  ?  &  ?  &  ?  &  $\ge 1411$  &    &    &    \\
121-128  &    &    &    &  ?  &  ?  &  ?  &  ?  &  ?  &  ?  &  ?  &    &    &    \\
129  &    &    &    &    &    &  ?  &  ?  &  ?  &  ?  &  ?  &  $\ge 1$  &    &    \\
130  &    &    &    &    &    &  ?  &  ?  &  ?  &  ?  &  ?  &  $\ge 4$  &    &    \\
131  &    &    &    &    &    &    &  ?  &  ?  &  ?  &  ?  &  $\ge 7$  &    &    \\
132-139  &    &    &    &    &    &    &  ?  &  ?  &  ?  &  ?  &  ?  &    &    \\
140  &    &    &    &    &    &    &    &    &  ?  &  ?  &  ?  &  1  &    \\
141  &    &    &    &    &    &    &    &    &  ?  &  ?  &  ?  &  0  &    \\
142-146  &    &    &    &    &    &    &    &    &  ?  &  ?  &  ?  &  ?  &    \\
147  &    &    &    &    &    &    &    &    &    &  ?  &  ?  &  $\ge 10$  &    \\
148  &    &    &    &    &    &    &    &    &    &  ?  &  ?  &  $\ge 28$  &    \\
149  &    &    &    &    &    &    &    &    &    &    &  ?  &  $\ge 39$  &    \\
150  &    &    &    &    &    &    &    &    &    &    &  ?  &  $\ge 27$  &    \\
151  &    &    &    &    &    &    &    &    &    &    &  ?  &  $\ge 19$  &    \\
152  &    &    &    &    &    &    &    &    &    &    &  ?  &  $\ge 11$  &    \\
153  &    &    &    &    &    &    &    &    &    &    &  ?  &  $\ge 8$  &    \\
154-155  &    &    &    &    &    &    &    &    &    &    &    &  ?  &    \\
156  &    &    &    &    &    &    &    &    &    &    &    &  ?  &  5  \\
157  &    &    &    &    &    &    &    &    &    &    &    &  ?  &  $\ge 6$ \\
158  &    &    &    &    &    &    &    &    &    &    &    &    &  $\ge 10$  \\
159  &    &    &    &    &    &    &    &    &    &    &    &    &  $\ge 6$  \\
160  &    &    &    &    &    &    &    &    &    &    &    &    &  $\ge 5$  \\
161  &    &    &    &    &    &    &    &    &    &    &    &    &  $\ge 2$  \\
162  &    &    &    &    &    &    &    &    &    &    &    &    &  $\ge 6$  \\
\hline
\end{tabular}
}

\caption{Number of $(3,J_{10};n,e)$-graphs, for $n \ge 24$.}

\label{table:graph_counts_K10}
\end{center}
\end{table}

\begin{table}[H]
\begin{center}
{\scriptsize
\begin{tabular}{| c | cccccc |}
\hline 
edges & \multicolumn{6}{c}{number of vertices $n$} \vline\\
  $e$ & 29 & 30 & 31 & 32 & 33 & 34\\
\hline 
58  &  5  &    &    &    &    &    \\
59  &  1364  &    &    &    &    &    \\
60-65  &  ?  &    &    &    &    &    \\
66  &  ?  &  5084  &    &    &    &    \\
67-72  &  ?  &  ?  &    &    &    &    \\
73  &  ?  &  ?  &  2657  &    &    &    \\
74-79  &  ?  &  ?  &  ?  &    &    &    \\
80  &  ?  &  ?  &  ?  &  4  &    &    \\
81  &  ?  &  ?  &  ?  &  6601  &    &    \\
82-89  &  ?  &  ?  &  ?  &  ?  &    &    \\
90  &  ?  &  ?  &  ?  &  ?  &  57099  &    \\
91-98  &  ?  &  ?  &  ?  &  ?  &  ?  &    \\
99  &  ?  &  ?  &  ?  &  ?  &  ?  &  $\ge 1$  \\
$\ge 100$  &  ?  &  ?  &  ?  &  ?  &  ?  &  ?  \\
\hline
\end{tabular}
}

\caption{Number of $(3,J_{11};n,e)$-graphs, for $29 \le n \le 34$.}
\label{table:graph_counts_K11}
\end{center}
\end{table}

\begin{table}[H]
\begin{center}
{\scriptsize
\begin{tabular}{| c | cccccc |}
\hline 
edges & \multicolumn{6}{c}{number of vertices $n$} \vline\\
  $e$ & 21 & 22 & 23 & 24 & 25 & 26\\
\hline 
45  &  1  &    &    &    &    &    \\
46  &  2  &    &    &    &    &    \\
47  &  61  &    &    &    &    &    \\
48  &  3743  &    &    &    &    &    \\
49  &  408410  &    &    &    &    &    \\
50-53  &  ?  &    &    &    &    &    \\
54  &  ?  &  2  &    &    &    &    \\
55  &  ?  &  299  &    &    &    &    \\
56  &  ?  &  20314  &    &    &    &    \\
57  &  ?  &  985296  &    &    &    &    \\
58  &  ?  &  23618486  &    &    &    &    \\
59-60  &  ?  &  ?  &    &    &    &    \\
63  &  ?  &  ?  &  9  &    &    &    \\
64  &  ?  &  ?  &  528  &    &    &    \\
65  &  ?  &  ?  &  24860  &    &    &    \\
66  &  ?  &  ?  &  566836  &    &    &    \\
67  &  ?  &  ?  &  5830123  &    &    &    \\
68-72  &  ?  &  ?  &  ?  &    &    &    \\
72  &  ?  &  ?  &  ?  &  104  &    &    \\
73  &  ?  &  ?  &  ?  &  1068  &    &    \\
74  &  ?  &  ?  &  ?  &  7913  &    &    \\
75  &  ?  &  ?  &  ?  &  31134  &    &    \\
76  &  ?  &  ?  &  ?  &  84634  &    &    \\
77  &  ?  &  ?  &  ?  &  160815  &    &    \\
78  &  ?  &  ?  &  ?  &  215365  &    &    \\
79  &  ?  &  ?  &  ?  &  207752  &    &    \\
80  &  ?  &  ?  &  ?  &  172746  &  18  &    \\
81  &  ?  &  ?  &  ?  &  142474  &  97  &    \\
82  &  ?  &  ?  &  ?  &  121641  &  333  &    \\
83  &  ?  &  ?  &  ?  &  107869  &  516  &    \\
84  &  ?  &  ?  &  ?  &  108001  &  416  &    \\
85  &    &  ?  &  ?  &  110938  &  158  &    \\
86  &    &  ?  &  ?  &  101090  &  30  &    \\
87  &    &  ?  &  ?  &  72665  &  5  &    \\
88  &    &  ?  &  ?  &  40935  &  1  &    \\
89  &    &    &  ?  &  17722  &  0  &    \\
90  &    &    &  ?  &  6262  &  0  &  3  \\
91  &    &    &  ?  &  1779  &  0  &  2  \\
92  &    &    &  ?  &  522  &  5  &  1  \\
93  &    &    &    &  129  &  16  &  0  \\
94  &    &    &    &  46  &  35  &  0  \\
95  &    &    &    &  8  &  34  &  0  \\
96  &    &    &    &  5  &  19  &  0  \\
97  &    &    &    &    &  6  &  0  \\
98  &    &    &    &    &  2  &  0  \\
99  &    &    &    &    &  1  &  0  \\
100  &    &    &    &    &  1  &  0  \\
101-103  &    &    &    &    &    &  0  \\
104  &    &    &    &    &    &  2  \\
\hline
$|\mathcal{R}(3,K_9-2K_2;n)|$  &  ?  &  ?  &  ?  &  1713617  &  1693  &  8  \\
\hline
\end{tabular}
}

\caption{Number of $(3,K_9-2K_2;n,e)$-graphs, for $n \ge 21$.}
\label{table:graph_counts_K9-2e}
\end{center}
\end{table}

\endgroup

\eject
\section*{Appendix 2: Testing correctness}
\label{section:appendix2}

Since most results obtained in this paper rely on computations,
it is very important that the correctness of our programs has
been thoroughly verified. Below we describe how we tested the correctness of our programs.

\bigskip
\noindent
{\bf Correctness}

\begin{itemize}
\item
For every $(3,J_k)$-graph which was output by our programs,
we verified that it does not contain a spanning subgraph of
$\overline{J_k}$ as induced subgraph by using an independent
program.

\item
Every Ramsey graph for $K_k$ is also a Ramsey graph for $J_{k+1}$.
Therefore, we verified that the complete lists of
$(3,J_{k+1};n,e)$-graphs which were generated by our programs
include all known $(3,k;n,e)$-graphs which we had found in~\cite{GoRa}.

\item
For every $(3,J_k;n,e(3,J_k,n))$-graph which was generated by our
programs, we verified that dropping any edge results in a graph
which contains a spanning subgraph of $\overline{J_k}$ as induced subgraph.

\item
For various $(3,J_k;n,\le e)$-graphs we added up to $f$ edges in all
possible ways to obtain $(3,J_k;n,\le e + f)$-graphs. For the cases
where we already had the complete set of $(3,J_k;n,\le e + f)$-graphs,
we verified that no new $(3,J_k;n,\le e + f)$-graphs were obtained.
We used this, amongst other cases, to verify that no new
$(3,J_{10};26,\le 55)$, $(3,J_{10};28,\le 70)$, $(3,J_{10};30,\le 87)$
or $(3,J_{11};32,\le 81)$-graphs were obtained.

\item
For various $(3,J_k;n,\le e + f)$-graphs we
dropped one edge in all
possible ways and verified that no new $(3,J_k;n,\le e + f - 1)$-graphs
were obtained.
We used this technique, amongst other cases, to verify that no new
$(3,J_{10};26,\le 54)$, $(3,J_{10};28,\le 69)$, $(3,J_{10};30,86)$
or $(3,J_{11};32,80)$-graphs were obtained.

\item
For various sets of $(3,J_{k+1};n,\le e)$-graphs we took each member
$G$ and constructed from it all $G_v$'s. We then
verified that this did not yield any new
$(3,J_k;n - deg(v) - 1,\le e - Z(v))$-graphs for the cases where
we have all such graphs.
We performed this test, amongst other cases, on the sets of
$(3,J_9;28,\le 70)$- and $(3,J_{10};32,\le 81)$-graphs.

\end{itemize}

Various sets of graphs can be obtained by both the maximum
triangle-free and neighborhood gluing extension method. Therefore,
as a test for the correctness of our implementations, we applied
both methods for the generation of several sets of graphs. We also
compared our results with known results. In each case, the results
were in complete agreement. More details are given below:

\begin{itemize}

\item
The sets of $(3,J_8;19,\le 38)$, $(3,J_8;20,\le 46)$,
$(3,J_8;21,\le 54)$ and $(3,J_9;27,\le 86)$-graphs were obtained
by both the maximum triangle-free method and the neighborhood gluing
extension method. The results were in complete agreement.

\item
The counts of all $(3,J_7)$-graphs are confirmed by~\cite{Fid}.

\item
The counts of all $(3,J_8;19,37)$, $(3,J_8;20,44)$, $(3,J_8;21,\le 52)$,
$(3,J_8;22,\le 60)$,\\
$(3,J_8;23,\le 71)$ and $(3,J_8;24,\le 81)$-graphs
are confirmed by~\cite{Ra1}.

\item
The counts of $(3,J_8;22,\le 65)$, $(3,J_8;23)$ and
$(3,J_8;24)$-graphs are confirmed by~\cite{MPR}.
\end{itemize}

Since our results are in complete agreement with previous results
and since all of our consistency tests passed, we believe that this
is strong evidence for the correctness of our implementations and results.

\bigskip
\noindent
{\bf Computation Time}

\medskip
We implemented the extension algorithms described
in Section~\ref{section:algorithms} in C.
Most computations were performed on a cluster
with Intel Xeon L5520 CPU's at 2.27 GHz, on which
a computational effort of one CPU year can be usually
completed in about 8 elapsed hours. The overall
computational effort which was required to improve the upper bounds of $R(3,J_k)$
is estimated to be about 40 CPU years. This includes
the time used by a variety of programs. The most CPU-intensive
task was the computation to determine all $(3,J_9;\ge 28)$-graphs
with the maximum triangle-free method. This took approximately 13 CPU years.
Also the computations to determine new lower bounds on $e(3,J_{11},n)$ took
relatively long. For example, it took nearly 5 CPU years using the neighborhood
gluing method to prove that 
$e(3,J_{11},39) \ge 151$.

The CPU time needed to complete the computations of Section~5 was negligible.